\documentclass[titlepage,draft,12pt]{article} 
\usepackage{amsfonts,latexsym,pstricks,pst-plot} 
\usepackage[a4paper]{geometry}

\date{}

\newcommand{\qed}{{\penalty 10000\mbox{$\quad\Box$}}}
\newcommand{\re}{\mathbb{R}}
\newcommand{\n}{\mathbb{N}}
\newcommand{\Gkn}{G^{(k)}_n}
\newcommand{\Gki}{G^{(k)}_\infty}

\newcommand{\LL}{L^2}
\newcommand{\PCD}{PC_{D}}
\newcommand{\PSD}{PS_{D}}
\newcommand{\PCn}{PC_{n}}
\newcommand{\PSDn}{PS_{D,n}}
\newcommand{\tsing}{T_{\mbox{{\scriptsize \textup{sing}}}}}
\newcommand{\tsingn}{T_{\mbox{{\scriptsize \textup{sing}}},n}}

\newtheorem{thm}{Theorem}[section]
\newtheorem{thmbibl}{Theorem}
\newtheorem{rmk}[thm]{Remark}
\newtheorem{prop}[thm]{Proposition}
\newtheorem{defn}[thm]{Definition}

\newtheorem{lemma}[thm]{Lemma}

\title{Slow time behavior of the semidiscrete Perona-Malik scheme in
dimension one}

\author{Maria Colombo\vspace{1ex}\\ 
{\normalsize Scuola Normale Superiore}\\
{\normalsize PISA (Italy)}\\
{\normalsize e-mail: \texttt{maria.colombo@sns.it}}
\and
Massimo Gobbino\vspace{1ex}\\ 
{\normalsize Universit\`a degli Studi di Pisa} \\
{\normalsize Dipartimento di Matematica Applicata ``Ulisse Dini''}\\ 
{\normalsize PISA (Italy)}\\  
{\normalsize e-mail: \texttt{m.gobbino@dma.unipi.it}}}

\begin{document}
\maketitle
\begin{abstract}
	We consider the long time behavior of the semidiscrete scheme for
	the Perona-Malik equation in dimension one.  We prove that
	approximated solutions converge, in a slow time scale, to
	solutions of a limit problem.  This limit problem evolves
	piecewise constant functions by moving their plateaus in the
	vertical direction according to a system of ordinary differential
	equations.
	
	Our convergence result is \emph{global-in-time}, and this forces
	us to face the collision of plateaus when the system singularizes.
	
	The proof is based on energy estimates and gradient-flow
	techniques, according to the general idea that ``the limit of the
	gradient-flows is the gradient-flow of the limit functional''.
	Our main innovations are a uniform H\"{o}lder estimate up to the
	first collision time included, a well preparation result with a
	careful analysis of what happens at discrete level during
	collisions, and renormalizing the functionals after each collision
	in order to have a nontrivial Gamma-limit for all times.

\vspace{1cm}

\noindent{\bf Mathematics Subject Classification 2000 (MSC2000):}
35K55, 35B40, 49M25.

\vspace{1cm} 

\noindent{\bf Key words:} Perona-Malik equation, semidiscrete scheme,
forward-backward parabolic equation, gradient-flow, maximal slope
curves, Gamma-convergence.
\end{abstract}

 
\section{Introduction}

The one dimensional Perona-Malik equation is the partial differential 
equation
\begin{equation}
	u_{t}=\left(\frac{u_{x}}{1+u_{x}^{2}}\right)_{x}=
	\frac{1-u_{x}^{2}}{(1+u_{x}^{2})^{2}}\, u_{xx}
	\quad\quad
	(x,t)\in(0,1)\times(0,+\infty),
	\label{eqn:PM}
\end{equation}
which is usually coupled with Neumann boundary conditions
\begin{equation}
	u_{x}(0,t)=u_{x}(1,t)=0
	\quad\quad
	\forall t>0,
	\label{eqn:nbc}
\end{equation}
and an initial condition
\begin{equation}
	u(x,0)=u_{0}(x)
	\quad\quad
	\forall x\in(0,1).
	\label{eqn:initial}
\end{equation}

Problem (\ref{eqn:PM}), (\ref{eqn:nbc}), (\ref{eqn:initial}) is the 
formal gradient-flow of the functional
\begin{equation}
	PM(u):=\frac{1}{2}\int_{0}^{1}\log\left(1+u_{x}^{2}\right)\,dx.
	\label{defn:PM-funct}
\end{equation}

The convex-concave behavior of the integrand in (\ref{defn:PM-funct}) 
makes (\ref{eqn:PM}) a parabolic equation of forward-backward type, with 
forward (or subcritical) regime in the region where $|u_{x}|<1$, and 
backward (or supercritical) regime in the region where $|u_{x}|>1$.

The analogous problem in two space dimensions was introduced by
P.~Perona and J.~Malik~\cite{PM} in the context of image denoising.
The rough idea is that small disturbances, corresponding to small
values of the gradient, are expected to be smoothed out by the
diffusion in forward regions.  On the contrary, sharp edges should be
enhanced by the backward nature of the equation in regions where the
gradient is large.

This intuition has actually been confirmed by numerical experiments.
There are some well known shortcomings, such as the staircasing effect
observed in supercritical regions, but nevertheless the method reveals
some stability, and in any case much more stability than expected from
a backward diffusion process.

Equation (\ref{eqn:PM}) is the prototype of all forward-backward 
parabolic equations such as $u_{t}=(\varphi'(u_{x}))_{x}$, where 
$\varphi$ is a nonconvex integrand. It is also strongly related to 
forward-backward parabolic equations of the form 
$u_{t}=(\phi(u))_{xx}$, where $\phi$ is a nonmonotone response 
function (indeed this is the equation solved by the derivative 
$u_{x}$ of solutions $u$ of (\ref{eqn:PM})). Such equations  
attracted a considerable attention in the last years because they are 
involved in several models, from phase transitions to population 
dynamic (see~\cite{tesei} and the references quoted therein).

A natural approach to an ill-posed problem is to approximate it by
more stable ones.  Following this idea, several authors proved well
posedness results for approximations of (\ref{eqn:PM}) obtained via
space discretization~\cite{GG} or convolution~\cite{CLMC}, time
delay~\cite{amann}, fractional derivatives~\cite{guidotti,gl}, fourth
order regularization~\cite{BF}, simplified nonlinearities~\cite{BNP}.
A satisfactory understanding of what happens as the suitable parameter
vanishes still seems to be out of reach.

Several papers reported numerical experiments on the Perona-Malik
equation in dimension one or two.  We refer in particular to
\cite{amann,BNPT1,E1,E2,guidotti,gl}.  All these experiments,
although obtained through different approximation methods, seem to
reveal some common \emph{qualitative} features.  In particular, the
evolution seems to happen in three different times scales.  We call
them ``fast time'', ``standard time'', and ``slow time'', according to
the terminology introduced in~\cite{BF}.
\begin{itemize}
	\item \emph{Fast time}.  In a time interval of order $o(1)$,
	solutions of approximated problems tend to develop microstructures
	in the concave region, with fast oscillations between very small
	and very large values of the derivative.  This is the staircasing
	effect, which causes an instantaneous drastic reduction of the
	energy in the backward regime.  From the variational point of view
	this is hardly surprising, due to the concavity of the integrand
	in that region.  More surprising is that this effect does not
	extend immediately to the forward regime, as it could be expected
	after remarking that the relaxation of (\ref{defn:PM-funct}) is
	trivially zero.
	
	Up to our knowledge, there is no rigorous treatment of this 
	phenomenon.  On the other hand, the existence of
	dense classes of smooth solutions of (\ref{eqn:PM})
	(see~\cite{GG-tams,GG-cpdes}) suggests that it is not reasonable
	to expect the staircasing effect for \emph{all} initial data with
	both subcritical and supercritical regions, but at most for
	``generic'' such data.  This remains a challenging open problem.

	\item \emph{Standard time}.  In a time interval of order $O(1)$,
	solutions of approximated problems evolve in order to reduce the
	energy in the convex region.  Rigorous results in this time scale
	are known only for the semidiscrete scheme in dimension one (in
	this paper semidiscrete means discrete with respect to space, and
	continuous with respect to time).  In~\cite{GG} a compactness
	result was proven, according to which solutions of approximated
	problems converge to something (as the size of the grid goes to
	zero), and all possible limits are classical solutions of
	(\ref{eqn:PM}) in the subcritical region of $u_{0}$.
	
	The characterization of such limits in supercritical regions
	remains an open problem, as well as any compactness result for
	different approximation methods or in more space dimensions.
 
	\item \emph{Slow time}.  After the second phase of the evolution,
	the energy has been reduced almost to zero, and the solution is
	close to a piecewise constant function.  This is consistent with
	the intuitive idea that piecewise constant functions are
	stationary points of $PM(u)$.  Since there is almost no
	energy left, the evolution slows down.
	
	Nevertheless, in a slower time scale the plateaus of this
	piecewise constant function tend to move in the vertical
	direction, with jump points which remain fixed in space.  The
	vertical dynamic is nontrivial because neighboring plateaus can
	collide, and actually do collide in a finite time.  After each
	collision at least one discontinuity point disappears, and the
	evolution proceeds as soon as the solution becomes a constant and
	there is nothing else to evolve.
\end{itemize}

The aim of this paper is a rigorous analysis of the slow time for the
semidiscrete scheme in dimension one (we refer to
section~\ref{sec:sds} for precise definitions).  The first three steps
in this direction were done by G.~Bellettini, M.~Novaga and M.~Paolini
in~\cite{BNP3}.  First of all, they identified the right time-scale,
which turns out to be of the same order of the inverse of the grid
size (namely $O(n)$ if the grid size is $1/n$).  Secondly, they
identified the system of ordinary differential equations describing
the evolution of the plateau heights in the limit problem.  Finally,
they proved that the rescaled solutions of the semidiscrete scheme
converge to the limit evolution described by that system in the
half-open interval $[0,\tsing)$, where $\tsing$ is the life span of 
the solution of the system.  

The proof of their convergence result is based on the construction of
suitable subsolutions and supersolutions, suggested by a formal
development of approximating solutions.  This method reveals some
drawbacks.  First of all, it requires some heavy computations, which
in~\cite{BNP3} are carried out at the expenses of choosing a
simplified form of the nonlinearity, a very special sequence of
initial data, and Dirichlet boundary conditions.  More important, it
seems quite hard to extend these arguments beyond $\tsing$, namely
when the interaction of plateaus makes the dynamic highly nontrivial.

In this paper we overcome these difficulties, and we prove a
\emph{global-in-time} convergence result (Theorem~\ref{thm:main}).  Of
course the limit problem is defined by restarting the evolution after
each collision according to the same rule applied to the new (smaller)
set of plateaus.  

Before restarting the evolution, it is however necessary to prove that
it can be extended up to $\tsing$ (included).  This fact
is quite intuitive, and indeed it has been implicitly mentioned (but
not proved) in~\cite{BNP3}, when the authors say that the system
singularizes due to collisions, and not to more strange phenomena.  On
the other hand, the possibility that more than two neighboring
plateaus collide in the same time makes this issue nontrivial.  We
overcome this difficulty by proving a 1/4-H\"{o}lder estimate up to
$\tsing$ (see Proposition~\ref{prop:u}, and Proposition~\ref{prop:hc}
for the corresponding estimate at discrete level).

Then we pass to our convergence result, inspired by the general
principle that ``the limit of gradient-flows is the gradient-flow of
the limit''.  Since approximating solutions are gradient-flows of
rescaled approximations of (\ref{defn:PM-funct}), it is reasonable to
expect the limit of the evolutions to be the gradient-flow of the
limit energy (in the sense of Gamma-convergence).  Unfortunately, if
we want the limit energy to be finite, we are forced to fix a priori
the number of discontinuities, and we are back to the interval
$[0,\tsing)$.

The first idea is therefore to renormalize the energy after each
collision.  If we add a constant (depending on the grid size) to each
approximated energy, then the approximated gradient-flows do not
change, but the limit energy can be different.  This allows to iterate
the convergence result after each collision, provided that we arrive
up to $\tsing$ (included) with approximating solutions which remain
``well prepared'', namely close enough in many senses to the
continuous limit.

To this end we develop two main tools.  The first one is a well
preparation result (Proposition~\ref{prop:wpl}).  When a discontinuity
disappears in the continuous limit, then in the corresponding interval
of the grid the discrete derivative of approximating solutions crosses
the critical threshold, switching from the concave region to the
convex one, and instantaneously the discrete solution becomes a well
prepared approximation of the new piecewise constant function with a
smaller set of plateaus.  The second tool
(Proposition~\ref{prop:fino-t1}) is a convergence result up to the
first collision time (included), which by itself improves the
convergence result of~\cite{BNP3}.  We prove it by rewriting both the
approximating problems, and the limit problem, in terms of
\emph{integral inequalities} instead of differential equations.  This
formulation, inspired by the theory of \emph{maximal slope curves}
(introduced in \cite{dgmt}, see~\cite{AGS} for a modern presentation),
happens to be much more stable when passing to the limit.

Our techniques work with general nonlinearities (we only need the
convex-concave behavior of the integrand), general sequences of
initial data (we do not even assume the boundedness of the energy),
and general boundary conditions (we work with Neumann boundary
conditions because this is the natural choice in applications, but the
same arguments apply to Dirichlet or periodic boundary conditions).

Of course several problems remain open.  Apart from the notorious
questions concerning fast time and standard time, it could be
interesting to prove similar results for the slow time in higher
dimension, or again in dimension one but with different approximation
methods.

A partial contribution in this direction is due to G.~Bellettini and
A.~Fusco~\cite{BF}.  They considered a fourth order regularization of
(\ref{eqn:PM}), corresponding to adding a vanishing second order term
to (\ref{defn:PM-funct}).  They identified the time-scale of slow
time, they computed the Gamma-limit of the rescaled energies, and they
conjectured that the limit problem is the gradient-flow of the limit
energy.  Unfortunately in that case this remains a conjecture, since
up to now no convergence result (even before collisions) is known.

The limit conjectured in~\cite{BF} evolves once again piecewise
constant functions, but the law of the vertical motion is in their
case different (the system of ordinary differential equations is
similar, but with different exponents).  This suggests two remarks.
On the one hand the existence of a slow time vertical motion is a
qualitative feature which is intrinsic in the nature of
(\ref{eqn:PM}).  On the other hand, what exactly happens in the slow
time from the quantitative point of view does depend on the
approximation method.

This paper is organized as follows.  In Section~\ref{sec:notations} we
introduce the rescaled semidiscrete scheme, the variational setting,
and the limit problem.  We also recall the previous results which are
needed throughout this paper.  In Section~\ref{sec:main} we state our
main results.  In Section~\ref{sec:tools} we present the basic tools
of our analysis.  In Section~\ref{sec:proofs} we collect all proofs.

\setcounter{equation}{0}
\section{Notation and definitions}\label{sec:notations}

\subsection{Functional spaces}

\subparagraph{\textmd{\emph{Continuous setting}}}

The more general ambient space we consider is $L^{2}((0,1))$,
shortened to $\LL$ when it is clear that we are working in the interval
$(0,1)$.  We write $\|u\|_{L^{p}((0,1))}$, or simply $\|u\|_{p}$, to
denote the $p$-norm ($p\in[0,+\infty]$) of a function $u$, and
$\langle u,v\rangle$ to denote the scalar product of the functions $u$
and $v$ in the appropriate $L^{2}$ space.

Let $D\subseteq(0,1)$ be a finite set, and let $k:=|D|$.  The elements
of $D$ divide $(0,1)$ into $(k+1)$ subintervals.

We call $\PCD$ the space of functions which are constant in each
subinterval, with the agreement that the constant values in any two
neighboring subintervals are different.  In other words, elements of
$\PCD$ are \emph{piecewise constant} functions with exactly $k$ jump
points located in the discontinuity set $D$.

We call $\PSD$ the space of functions which are Lipschitz continuous
in each subinterval, with Lipschitz constant less than or equal to 1,
with the agreement that for each $d\in D$ the limit as $x\to d^{-}$ is
different from the limit as $x\to d^{+}$.  In other words, elements of
$\PSD$ are \emph{piecewise subcritical} functions with exactly $k$
jump points located in the discontinuity set $D$.  For every
$u\in\PSD$, and every $d\in D$, the \emph{jump height} of $u$ in $d$
is defined as
\begin{equation}
	J_{d}(u):=\lim_{x\to d^{+}}u(x)-\lim_{x\to d^{-}}u(x).
	\label{defn:Jd}
\end{equation}

It is easy to see that $\PCD\subseteq\PSD\subseteq\LL$.  Every element
of $\PCD$ is uniquely determined by the heights of its $(k+1)$
plateaus. This correspondence defines an isometry between $\PCD$ and 
an open subset of a Euclidean space of dimension $(k+1)$. When 
needed, we assume that elements of $\PCD$ and $\PSD$ are defined in 
the jump points in such a way that they are right-continuous.

\subparagraph{\textmd{\emph{Discrete setting}}}

Given a positive integer $n$, we divide $[0,1]$ into $n$ intervals of
length $1/n$, and we consider the space $\PCn$ of all functions which
are constant in each subinterval (in this case constants in
neighboring subintervals may be equal).  The space $\PCn$, when
endowed with the $L^{2}$-norm inherited as a subset of $\LL$, becomes a
Euclidean space isomorphic to $\re^{n}$.  Since elements of $\PCn$ are
thought as $\LL$ functions, it is not so essential to define them also
in points of the form $i/n$ (with $i=0,1,\ldots,n$).  In any case, when
needed we assume that the value in any of these points is the same as
in the interval on its right (on its left in the case $i=n$).

Given $u\in\PCn$, the \emph{discrete derivatives} $D^{1/n}u$ and
$D^{-1/n}u$ are defined as the incremental quotients
$$D^{\pm 1/n}u(x):=\frac{u(x\pm 1/n)-u(x)}{\pm 1/n}
\quad\quad
\forall x\in[0,1],$$
with the agreement that $u$ has been extended previously to the whole
real line (or at least to a neighborhood of $[0,1]$ of width $1/n$) by
setting $u(x)=u(0)$ for every $x\leq 0$, and $u(x)=u(1)$ for every
$x\geq 1$.  

Given a finite set $D\subseteq(0,1)$ with $k$ elements, we set
\begin{equation}
	D_{n}:=\bigcup_{d\in D}\left[\frac{\lceil nd\rceil-1}{n},
	\frac{\lceil nd\rceil}{n}\right)\subseteq[0,1].
	\label{defn:Dn}
\end{equation}

In other words, $D_{n}$ is the union of all subintervals which
intersect $D$ (when $d$ is of the form $i/n$ we take the subinterval
on its left).  It is easy to see that, when $n$ is large enough,
$D_{n}$ is the union of $k$ disjoint intervals, and $(0,1)\setminus
D_{n}$ has $(k+1)$ connected components. Since we are interested in 
passing to the limit as $n\to +\infty$, we can always work under this 
assumption.

We call $\PSDn$ the set of all functions $u\in\PCn$ such that
\begin{equation}
	\left|D^{1/n}u(x)\right|\leq 1\Longleftrightarrow x\in[0,1]\setminus
	D_{n}.
	\label{defn:PSDn}
\end{equation}

The space $\PSDn$ is obviously the discrete counterpart of $\PSD$.  In
analogy with (\ref{defn:Jd}), the \emph{discrete jump height} of a
function $u\in\PSDn$ in a point $d\in D$ is defined as
\begin{equation}
	J_{d,n}(u):=u\left(\frac{\lceil nd\rceil}{n}\right)-
	u\left(\frac{\lceil nd\rceil-1}{n}\right).
	\label{defn:Jdn}
\end{equation}

This is equivalent to say that $J_{d,n}(u):=u(x+1/n)-u(x)$, where $x$
is any point of the subinterval containing $d$, or of the subinterval
on the left of $d$ if $d=i/n$ for some $i=1,\ldots,n-1$.  We point out
that (\ref{defn:Jdn}) makes sense for every $d\in[0,1]$, and not only
for $d\in D$.  Of course we have that $|u(x+1/n)-u(x)|>1/n$ if and
only if $x\in D_{n}$.

The \emph{subcritical incremental quotient} of a function $u\in\PSDn$
is defined as
\begin{equation}
	SQ_{n}(u):=\left\|D^{1/n}u(x)\right\|_{L^{\infty} ((0,1)\setminus
	D_{n})}.
	\label{defn:SQn}
\end{equation}

Due to (\ref{defn:PSDn}) we have that $0\leq SQ_{n}(u)\leq 1$ for
every $u\in\PSDn$.  The subcritical incremental quotient is the
discrete counterpart of the Lipschitz constant in the intervals
between discontinuities.

\subsection{Functionals}

A discrete approximation of (\ref{defn:PM-funct}) is obtained by
replacing the derivative with discrete derivatives.  Thus we introduce
the functionals $PM_{n}:\PCn\to\re$ defined by
\begin{equation}
	PM_{n}(u)=\frac{1}{2}\int_{0}^{1}\log\left(1+
	|D^{1/n}u(x)|^{2}\right)dx
	\quad\quad
	\forall u\in\PCn.
	\label{defn:PMn}
\end{equation}

The time rescaling due to the ``slow time'' leads us to consider also the
functionals $n PM_{n}(u)$.  More generally, for each nonnegative integer
$k$ we consider the sequence of functionals $\Gkn:\PCn\to\re$ defined
by
\begin{equation}
	\Gkn(u):=\frac{n}{2}\int_{0}^{1}
	\log\left(1+|D^{1/n}u(x)|^{2}\right)dx-k\log n \quad\quad
	\forall u\in\PCn.
	\label{defn:gkn}
\end{equation}

It is clear that $\Gkn(u)$ (which sometimes we call \emph{$k$-energy}
of $u$) coincides with $n PM_{n}(u)$ up to an additive constant, and
in particular these functionals have the same gradient, hence also the
same gradient-flow.

Computing the gradient is a simple exercise in 
finite dimension. It turns out that
\begin{equation}
	\nabla\Gkn(u)=n\nabla PM_{n}(u)=-n\,D^{-1/n}\left[
	\frac{D^{1/n}u}{1+|D^{1/n}u|^{2}}\right]
	\quad\quad
	\forall u\in\PCn.
	\label{defn:grad-gkn}
\end{equation}

The ``second order discrete operator'' in the right-hand side of
(\ref{defn:grad-gkn}) needs some interpretation in the two extremal
subintervals, where it requires to compute values of $u$ outside
$[0,1]$.  Once again this is done after extending $u$ by setting
$u(x)=u(0)$ for every $x\leq 0$, and $u(x)=u(1)$ for every $x\geq 1$.
This agreement is the discrete counterpart of the Neumann boundary
condition, which in this sense is now included in the right-hand side
of (\ref{defn:grad-gkn}).

Finally, for every finite set $D\subseteq(0,1)$ with $|D|=:k$, and 
every $u\in\PSD$, we set
\begin{equation}
	\Gki(u):=\sum_{d\in D}\log|J_{d}(u)|.
	\label{defn:gki}
\end{equation}

The following result justifies the notation used for $\Gki$, and shows
that the sequence $\Gkn(u)$ has a less trivial limit as $n\to +\infty$
with respect to the sequence $n PM_{n}(u)$.  A proof of this result
can be found in~\cite{BNPT2}, or simply deduced from the theory of
convex-concave integrands developed in \cite{mora}.

	\begin{thmbibl}[Gamma-limit of discrete
	functionals]\label{thmbibl:gamma-limit} 
	Let $k$ be a nonnegative
	integer.  Let us extend $PM_{n}$ and $\Gkn$ by setting them equal
	to $+\infty$ for every $u\in\LL\setminus\PCn$.
	
	Then we have that (all $\Gamma$-limits are intended with respect 
	to $\LL$-metric)
	$$\Gamma\!-\!\!\!\lim_{n\to +\infty}nPM_{n}(u)=\left\{
	\begin{array}{ll}
		0 & \mbox{if $u$ is constant,}  \\
		+\infty & \mbox{otherwise,}
	\end{array}
	\right.$$
	and 
	$$\Gamma\!-\!\!\!\lim_{n\to +\infty}\Gkn(u)=\left\{
	\begin{array}{ll}
		-\infty & \mbox{if $u\in\PCD$ for some $D\subseteq(0,1)$ 
		with $|D|<k$,}  \\
		\noalign{\vspace{1ex}}
		\Gki(u) & \mbox{if $u\in\PCD$ for some $D\subseteq(0,1)$ 
		with $|D|=k$,}  \\
		\noalign{\vspace{1ex}}
		+\infty & \mbox{otherwise}.
	\end{array}
	\right.$$
\end{thmbibl}

\subsection{The semidiscrete scheme}\label{sec:sds}

The semidiscrete scheme for the one dimensional Perona-Malik equation 
is the gradient-flow of (\ref{defn:PMn}). This leads to the problem
\begin{equation}
	v_{n}'(t)=-\nabla PM_{n}(v_{n}(t))
	\quad\quad
	\forall t\geq 0,
	\label{pbm:semidiscrete}
\end{equation}
\begin{equation}
	v_{n}(0)=u_{0n},
	\label{pbm:data}
\end{equation}
where $\{u_{0n}\}$ is a suitable sequence of initial conditions with 
$u_{0n}\in\PCn$ for every $n\geq 1$. The behavior of $v_{n}(t)$ as 
$n\to +\infty$ is the subject of the ``standard time'' theory.

In the ``slow time'' theory we speed up the evolution by considering 
the sequence $u_{n}(t):=v_{n}(nt)$. It is very simple to show that 
this sequence solves the rescaled problems
\begin{equation}
	u_{n}'(t)=-n\nabla PM_{n}(u_{n}(t))
	\quad\quad
	\forall t\geq 0,
	\label{pbm:main-eq}
\end{equation}
\begin{equation}
	u_{n}(0)=u_{0n}.
	\label{pbm:main-data}
\end{equation}

Thanks to (\ref{defn:grad-gkn}), the differential equation
(\ref{pbm:main-eq}) is equivalent to
\begin{equation}
	u_{n}'(t)=-\nabla\Gkn(u_{n}(t))
	\quad\quad
	\forall t\geq 0.
	\label{pbm:main-eq-gkn}
\end{equation}

All these problems admit a unique solution defined for every $t\geq
0$.  Indeed $\PCn$ is a finite dimensional vector space, the
functionals we consider are of class $C^{\infty}$, and their gradient
(\ref{defn:grad-gkn}) is globally Lipschitz continuous.  It is
worthwhile to notice that formula (\ref{defn:grad-gkn}) makes
(\ref{pbm:semidiscrete}) the discrete counterpart of (\ref{eqn:PM}).

Throughout this paper we consider $u_{n}$ both as a function
$u_{n}(t)$ of the time variable with values in $\LL$, and as a function
$u_{n}(x,t)$ of $(x,t)$ with real values.  

The following properties of $u_{n}$ are used several times.  The proof
can be deduced from the corresponding properties of $v_{n}$ stated in
\cite{GG}.

\begin{thmbibl}[Properties of approximating 
	solutions]\label{thmbibl:semidiscrete}
	Let $n$ be a fixed positive integer, let $u_{0n}\in\PCn$, and let
	$u_{n}:[0,+\infty)\to\LL$ be the solution of problem
	(\ref{pbm:main-eq}), (\ref{pbm:main-data}).
	
	Then the following properties hold true.
	\begin{enumerate}
		\renewcommand{\labelenumi}{(\arabic{enumi})}
		\item  \emph{(Regularity)} We have that 
		$$u_{n}\in C^{\infty}\left([0,+\infty);\PCn\right)\subseteq
		C^{\infty}\left([0,+\infty);\LL\right).$$
	
		\item\label{stat:holder-1/2} \emph{(Standard gradient-flow
		estimate)} Let $k$ be any nonnegative integer.  Then the
		function $t\to\Gkn(u_{n}(t))$ is nonincreasing, and for every
		$0\leq s\leq t$ we have that
		$$\|u_{n}(t)-u_{n}(s)\|_{2}\leq\left\{
		\Gkn(u_{n}(s))-\Gkn(u_{n}(t))\right\}^{1/2}|t-s|^{1/2}.$$
	
		\item \emph{($L^{p}$ estimate)} The function
		$t\to\|u_{n}(x,t)\|_{L^{p}((0,1))}$ is nonincreasing for every
		$p\geq 1$ (including $p=\infty$).
	
		\item \emph{(Total variation estimate)}
		The function $t\to\|D^{1/n}u_{n}(x,t)\|_{L^{1}((0,1))}$ is
		nonincreasing.  The same is not necessarily true for $p$-norms
		with $p> 1$.
	
		\item\label{stat:asymptotic} \emph{(Asymptotic behavior)} The
		function $u_{n}(t)$ tends, as $t\to +\infty$, to the constant
		function equal to the average of $u_{0n}$.
	
		\item\label{stat:monotonicity} \emph{(Monotonicity and
		extinction of supercritical regions)} Supercritical regions
		are nonincreasing (as set valued maps), and they disappear
		after a finite time.
		
		In other words, if $u_{0n}\in\PSDn$ for some finite set
		$D\subseteq(0,1)$, then there exist $j\in\n$, and a finite
		sequence of times
		$$0=T_0<T_1<\ldots <T_{j}<T_{j+1}=+\infty,$$
		and a finite sequence of subsets
		$$D=D(0)\supset D(1)\supset\ldots\supset D(j)=\emptyset$$
		(with strict inclusions) such that
		$$u_n(t)\in PS_{D(i), n}
		\quad\quad
		\forall t\in [T_{i}, T_{i+1}),\ 
		\forall i\in \{0,1,\ldots,j\}.$$
	\end{enumerate}
\end{thmbibl}

The fact that supercritical regions disappear in a finite time (which
of course depends on $n$) has probably never been stated in the
literature, but it is a simple consequence of
statement~(\ref{stat:asymptotic}).  In turn,
statement~(\ref{stat:asymptotic}) follows from three general facts:
the average of $u_{n}(t)$ is invariant during the evolution, the limit
of a gradient-flow is a steady state solution, and
(\ref{defn:grad-gkn}) is zero if and only if $u$ is constant.

\subsection{The limit problem}\label{sec:limit-pbm}

Let $k$ be a nonnegative integer, and let $D\subseteq(0,1)$ be a 
finite set with $|D|=k$. Given an initial condition $u_{0}\in\PCD$, 
we define an evolution $u(t)$ starting from $u_{0}$ according to the 
following algorithm.

If $k=0$, then the initial datum $u_{0}$ is constant, and we define
$u(t)$ as the stationary solution $u(t)\equiv u_{0}$.

If $k>0$, let $D=\{d_{1},\ldots,d_{k}\}$ with
$0<d_{1}<\ldots<d_{k}<1$, let $d_{0}:=0$, $d_{k+1}:=1$, and let
$a_{0i}$ (with $i=0,\ldots,k$) denote the constant value of $u_{0}$ in
the interval $(d_{i},d_{i+1})$.  Let
$(a_{0}(t),a_{1}(t),\ldots,a_{k}(t))$ be the (unique) solution of the
system of $(k+1)$ ordinary differential equations
	\begin{eqnarray*}
		a_{0}'(t) & = & \frac{1}{d_{1}-d_{0}}\cdot
		\frac{1}{a_{1}(t)-a_{0}(t)},  \\
		a_{i}'(t) & = & \frac{1}{d_{i+1}-d_{i}}
		\left(\frac{1}{a_{i+1}(t)-a_{i}(t)}-  
		\frac{1}{a_{i}(t)-a_{i-1}(t)}\right)
		\quad\quad
		i=1,\ldots,k-1,\\
		a_{k}'(t) & = & -\frac{1}{d_{k+1}-d_{k}}\cdot
		\frac{1}{a_{k}(t)-a_{k-1}(t)},
	\end{eqnarray*}
with initial conditions $a_{i}(0)=a_{0i}$ for every $i=0,\ldots,k$. 

Let $u(t)$ be the piecewise constant function whose value in
$(d_{i},d_{i+1})$ is $a_{i}(t)$, defined as soon as the solution of 
the system exists. To this end, we have the following result (the 
proof is given is section~\ref{sec:proof-u}).

\begin{prop}\label{prop:u}
	Let $k$, $D$, $u_{0}$, $u(t)$ be as above.  Then we have the
	following conclusions.
	\begin{enumerate}
		\renewcommand{\labelenumi}{(\arabic{enumi})}
		\item  \emph{(Local but not global existence)} The system of 
		ordinary differential equations has a local solution defined 
		on a maximal interval $[0,\tsing)$ with
		$\tsing\in(0,+\infty)$.
	
		\item  \emph{($L^{\infty}$ estimate)} We have that
		\begin{equation}
			\|u(t)\|_{\infty}\leq\|u_{0}\|_{\infty} 
			\quad\quad 
			\forall t\in [0,\tsing).
			\label{th:infty-u}
		\end{equation}
	
		\item\label{stat:u-hc} \emph{(1/4-H\"{o}lder continuity up to
		collision)} For every $(s,t)\in [0,\tsing)^{2}$ we have that
		\begin{equation}
			\|u(t)-u(s)\|_{2}\leq (3k)^{3/4}\exp\left(
			\frac{1}{2k}\Gki(u_{0})\right)
			|t-s|^{1/4}.
			\label{th:h-1/4-u}
		\end{equation}
	\end{enumerate}
\end{prop}
Thanks to the H\"{o}lder continuity (\ref{th:h-1/4-u}), we can define
$u(t)$ up to $\tsing$ (included).  This extension fulfils
(\ref{th:h-1/4-u}) and (\ref{th:infty-u}) in the closed interval.

Moreover, $u(t)\in PC_{D'}$ for some $D' \subseteq D$.  If $D'=D$ we
can continue the solution of the system of ordinary differential
equations beyond $\tsing$, but this contradicts the maximality of
$\tsing$.  Thus $D'$ is strictly contained in $D$, which means that at
time $t=\tsing$ we have a collision between at least two adjacent
plateaus.
At this point we restart the construction of $u(t)$ from $u(\tsing)$,
which has a smaller set of jump points.

This procedure defines a function $u\in
C^{0}\left([0,+\infty);\LL\right)$ with $u(0)=u_{0}$.  For each $t\geq
0$ the function $u(t)$ is piecewise constant in the space variable,
and its discontinuity set is contained in the discontinuity set of
$u_{0}$.  There is a finite set of singular times when two or more
adjacent plateaus collide, hence one or more discontinuities disappear.
After each collision, the involved plateaus remain attached forever,
and the evolution goes on according to the same rule applied to the
new set of plateaus.  After the last collision $u(t)$ becomes
constant, and it does not move any more.  This constant is actually
the average of $u_{0}$ (indeed the average of $u(t)$ is invariant
during the evolution).  The function $u$ is of class $C^{\infty}$ with
respect to the time variable outside the finite set of collision
times, and uniformly continuous as a function from $[0,+\infty)$ to
$\LL$.

\setcounter{equation}{0}
\section{Main results}\label{sec:main}

The main result of this paper is the following convergence result.

\begin{thm}[Global-in-time convergence]\label{thm:main}
	Let $D'\subseteq D\subseteq(0,1)$ be two finite sets.  Let
	$u_{0}\in PC_{D'}$, and let $\{u_{0n}\}$ be a sequence such that
	\begin{eqnarray}
		 & u_{0n}\in\PSDn 
		 \quad\quad 
		 \forall n\geq 1, & 
		\label{hp:u0n-in-psdn}  \\
		\noalign{\vspace{1ex}}
		 & u_{0n}\to u_{0}\quad\mbox{in }\LL((0,1)). & 
		\label{hp:u0n-to-u0}  
	\end{eqnarray}
	
	For every $n\geq 1$, let $u_{n}(t)$ be the solution of the
	approximating problem (\ref{pbm:main-eq}), (\ref{pbm:main-data}).
	Let $u(t)$ be the solution of the limit problem with initial
	condition $u_{0}$, as defined in section~\ref{sec:limit-pbm}.
	
	Then we have the following conclusions.
	\begin{enumerate}
		\renewcommand{\labelenumi}{(\arabic{enumi})} \item
		\emph{(Global-in-time $\LL$-convergence)} We have that
		$u_{n}(t)\to u(t)$ in
		$C^{0}\left([0,+\infty);\LL\right)$, namely
		\begin{equation}
			\lim_{n\to +\infty}\sup_{t\geq 0}
			\|u_{n}(t)-u(t)\|_{\LL((0,1))}=0.
			\label{th:main}
		\end{equation}
	
		\item \emph{(Global-in-time ``uniform'' convergence)} For every
		$t\geq 0$, let $D(t)$ be the discontinuity set of $u(t)$, and
		let $D_{n}(t)$ be the union of all subintervals containing
		elements of $D(t)$, defined according to (\ref{defn:Dn}) with
		$D(t)$ instead of $D$.  Let us set
		\begin{equation}
			\mathcal{K}_{n}:=\{(x,t)\in[0,1]\times[0,+\infty): 
			x\not\in D_{n}(t)\}.
			\label{defn:Kn}
		\end{equation}
		
		Then we have that
		\begin{equation}
			\lim_{n\to +\infty}
			\|u_{n}(x,t)-u(x,t)\|_{L^{\infty}(\mathcal{K}_{n})}=0.
			\label{th:main-unif}
		\end{equation}
	\end{enumerate}
\end{thm}

We considered a piecewise constant initial datum $u_{0}$ because this
is the natural space where the Gamma-limit of the renormalized
functionals is finite.  On the other hand, we emphasize that the
approximating sequence $\{u_{0n}\}$ is quite general.  In particular,
we did not assume that it is a recovery sequence, or that its energy
is bounded, and the set of discrete jump points of $u_{0n}$ can be
larger than the set of jump points of $u_{0}$.  What is essential is
that all discrete jump points of approximating functions are contained
in a \emph{fixed} finite set $D$ (if not, there are counterexamples
even to local-in-time convergence).

Another reason for looking at piecewise constant data is that they are
expected to be the limit as $t\to +\infty$ of evolutions in ``standard
time''.  Up to our knowledge, this has been proved rigorously only for
``generic'' piecewise subcritical data (see~\cite{BNPT1,GG}).

Whenever the ``standard time'' evolution admits a piecewise constant
limit, we can start our ``slow time'' analysis from that limit.  We
state this idea formally in the next result.  We point out that in
this case we do not assume that the sequence $\{u_{0n}\}$ of initial
data has a piecewise constant or piecewise subcritical limit (actually
we do not even assume that it has a limit).

\begin{thm}[Convergence for more general initial data]\label{cor:main}
	Let $D\subseteq(0,1)$ be a finite set, and let $\{u_{0n}\}$ be a
	sequence satisfying (\ref{hp:u0n-in-psdn}).
	
	For every $n\geq 1$, let $v_{n}(t)$ be the solution of problem
	(\ref{pbm:semidiscrete}), (\ref{pbm:data}) (no time rescaling),
	and let $u_{n}(t)=v_{n}(nt)$ be the solution of the rescaled
	problem (\ref{pbm:main-eq}), (\ref{pbm:main-data}).  Let us assume
	that there exist $S\geq 0$, $v\in C^{0}([S,+\infty),\LL)$,
	$D'\subseteq D$, and $v_{\infty}\in PC_{D'}$ such that
	\begin{eqnarray}
		 & \displaystyle{\lim_{n\to +\infty}v_{n}(t)=v(t)}
		\quad\quad
		\forall t\geq S, & 
		\label{hp:cormain-lim-vn}  \\
		 & \displaystyle{\lim_{t\to +\infty}v(t)=v_{\infty},} & 
		\label{hp:cormain-lim-v}
	\end{eqnarray}
	where both limits are intended in $\LL$.  Let $u(t)$ be the
	solution of the limit problem, defined as in
	section~\ref{sec:limit-pbm}, with initial condition
	$u(0)=v_{\infty}$.
	
	Then for every $T>0$ we have that
	$u_{n}(t)\to u(t)$ in $C^{0}\left([T,+\infty);\LL\right)$, 
	namely
	\begin{equation}
		\lim_{n\to +\infty}\sup_{t\geq T}
		\|u_{n}(t)-u(t)\|_{\LL((0,1))}=0.
		\label{th:main-cor}
	\end{equation}
	
	Moreover, if $\mathcal{K}_{n}$ is defined as in (\ref{defn:Kn}),
	we have that 
	\begin{equation}
		\lim_{n\to +\infty} 
		\|u_{n}(x,t)-u(x,t)\|_{L^{\infty}\left(\mathcal{K}_{n}
		\cap([0,1]\times[T,+\infty))\right)}=0.	
		\label{th:main-unif-cor}
	\end{equation}
\end{thm}

We conclude by mentioning a possible extension of our results.

\begin{rmk}
	\begin{em}
		For the sake of simplicity, we devoted this paper to the model
		case of the Perona-Malik equation, in which the integrand is
		$\varphi(\sigma):=2^{-1}\log(1+\sigma^{2})$.  Similar
		techniques apply to larger classes of convex-concave
		integrands, for example the case where
		$\varphi(\sigma):=\alpha^{-1}(1+\sigma^{2})^{\alpha/2}$ for
		some $\alpha\in(0,1)$.
		In this case the ``slow time'' is of order $O(n^{1-\alpha})$, 
		the limit energy is
		\begin{equation}
			G_{\alpha}(u):= \sum_{d\in D}|J_{d}(u)|^{\alpha},
			\label{defn:g-alpha}
		\end{equation}
		and the system of ordinary
		differential equations governing the evolution of the 
		plateau heights is
		$$a_{i}'(t) =
		\frac{1}{d_{i+1}-d_{i}}
		\left(
		\frac{a_{i+1}(t)-a_{i}(t)}{\left|a_{i+1}(t)-a_{i}(t)\right|^{2-\alpha}}-
		\frac{a_{i}(t)-a_{i-1}(t)}{\left|a_{i}(t)-a_{i-1}(t)\right|^{2-\alpha}}\right),$$
		suitably modified for $i=0$ and $i=k$.  
		
		There are, however, some remarkable differences.  On the one
		hand, this situation is simpler because the limit
		energy~(\ref{defn:g-alpha}) is bounded from below, hence there
		is no need to renormalize it after each collision.  On the
		other hand, the limit energy can be finite even if $u$ has
		infinitely many jump points.  
	\end{em}
\end{rmk}

\setcounter{equation}{0}
\section{Fundamental tools}\label{sec:tools}

In this section we state the main ingredients needed in the proof of
our main results.

The first one is a well preparation result, which plays its role at
the beginning of the evolution and during each collision.  In input we
have a sequence satisfying (\ref{hp:u0n-in-psdn}) and
(\ref{hp:u0n-to-u0}) as in the assumptions of Theorem~\ref{thm:main}.
In particular, some of the jump points might disappear in the limit
(this happens if and only if $D'$ is strictly contained in $D$), and
there is no information on the $k$-energy or the $k'$-energy of the
sequence (where $k:=|D|$ and $k':=|D'|$).  We prove the existence of a
sequence of times $S_{n}\to 0$ such that $u_{n}(S_{n})$ is a ``well
prepared'' sequence, namely it still converges to $u_{0}$, all its
elements lie in the corresponding space $PS_{D',n}$, and their
$k'$-energies converge to the $k'$-energy of $u_{0}$.

\begin{prop}[Well preparation]\label{prop:wpl}
	Let $D$, $D'$, $\{u_{0n}\}$, $u_{0}$, $u_{n}(t)$ be as in
	Theorem~\ref{thm:main}, and let $k':=|D'|$.
	
	Then there exists a sequence $S_{n} \to 0$ of positive real numbers 
	such that
	\begin{eqnarray}
		 & u_{n}(S_{n})\in PS_{D',n}
		 \quad
		\mbox{for every $n$ large enough}, & 
		\label{th:wpl-dje}  \\
		\noalign{\vspace{1ex}}
		 & \displaystyle{\lim_{n\to 
		 +\infty}G^{(k')}_{n}(u_{n}(S_{n}))=
		 G^{(k')}_{\infty}(u_{0}),} & 
		\label{th:wpl-wp}  \\
		\noalign{\vspace{0.3ex}}
		 & \displaystyle{\lim_{n\to +\infty}\;\max_{t\in[0,S_{n}]}
		\|u_{n}(t)-u_{0}\|_{2}= 0.} & 
		\label{th:wpl-sup}
	\end{eqnarray}
\end{prop}

The second tool is a convergence result up to the first jump
extinction.  It plays its role in the time intervals between
collisions.  Now in input we have a ``well prepared'' sequence, or at
least a sequence of initial data with bounded energy, and without
vanishing jump points.  We prove some sort of uniform convergence on
an increasing sequence of time intervals (depending on $n$).  As $n\to
+\infty$ these intervals invade the whole time interval up to the
first collision.

\begin{prop}[Convergence up to first jump extinction]\label{prop:fino-t1}
	Let $k$ be a positive integer, and let $D\subseteq(0,1)$ be a
	finite set with $|D|=k$.  Let $u_{0}\in\PCD$, and let $\{u_{0n}\}$
	be a sequence satisfying (\ref{hp:u0n-in-psdn}),
	(\ref{hp:u0n-to-u0}), and 
	\begin{equation}
		\sup_{n\geq 1}\Gkn(u_{0n})<+\infty.
		\label{hp:main-gkn}
	\end{equation}
	
	For every $n\geq 1$, let $u_{n}(t)$ be the solution of the
	approximating problem (\ref{pbm:main-eq}), (\ref{pbm:main-data}),
	and let 
	\begin{equation}
		\tsingn:=\sup\left\{t\geq 0:u_{n}(t)\in\PSDn\right\}
		\label{defn:tsingn}
	\end{equation}
	be the first time when a discrete jump disappears (it is the time
	$T_{1}$ in statement~(\ref{stat:monotonicity}) of
	Theorem~\ref{thmbibl:semidiscrete}).  Let $u(t)$ and $\tsing$ be
	defined as in section~\ref{sec:limit-pbm}.
	
	Then there exists a sequence $\{T_{n}\}$ of real numbers such that
	\begin{eqnarray}
		 & 0<T_{n}<\tsingn
		 \quad\quad
		 \forall n\geq 1,& 
		\label{th:fino-t1-tn}  \\
		\noalign{\vspace{1ex}}
		 & \displaystyle{\lim_{n\to +\infty}T_{n}=\tsing,} & 
		\label{th:fino-t1-lim}  \\
		 & \displaystyle{\lim_{n\to +\infty}\,
		 \max_{t\in[0,T_{n}]}\|u_{n}(t)-u(t)\|_{2}=0.} & 
		\label{th:fino-t1-sup}
	\end{eqnarray}
	
	Moreover we have that
	\begin{equation}
		\lim_{n\to +\infty}\Gkn(u_{n}(t))=\Gki(u(t))
		\quad\quad
		\forall t\in(0,\tsing).
		\label{th:t1-gkn<}
	\end{equation}
\end{prop}

Finally, we present a qualitative property of approximating solutions
which could be interesting in itself.  It is the discrete analog of
statement~(\ref{stat:u-hc}) of Proposition~\ref{prop:u}.  We point out
that the standard gradient-flow estimates (statement
(\ref{stat:holder-1/2}) of Theorem~\ref{thmbibl:semidiscrete}) control
the 1/2-H\"{o}lder constant of $u_{n}(t)$ in terms of the descent of
the energy, but such estimates are useless if the energy is not
bounded from below independently on $n$, and this is exactly what
happens in this model when $t$ approaches a collision time.

The following 1/4-H\"older estimates overcome this difficulty.

\begin{prop}[Uniform 1/4-H\"older continuity]\label{prop:hc} 

	Let $k$, $D$, $\{u_{0n}\}$, $u_{n}(t)$, $\tsingn$ be as in
	Proposition~\ref{prop:fino-t1}.
	
	Then for every $n\geq 1$ and every $(s,t)\in [0,\tsingn]^{2}$ we
	have that
	$$\|u_n(t)-u_n(s)\|_{2}\leq (3k)^{3/4}\exp\left(
	\frac{1}{2k}\Gkn(u_{0n})\right)
	|t-s|^{1/4}.$$
\end{prop}

\setcounter{equation}{0}
\section{Proofs}\label{sec:proofs}

\subsection{Basic estimates}

In this section we collect some general facts, which are going to be
used several times in the proof of our main results and basic tools.
The first one concerns ``double limits'' (we omit the simple proof).

\begin{lemma}[Double index sequence]\label{lemma:amn}
	Let $\{A_{m,n}\}$ (with $(m,n)\in\n^{2}$) be a double indexed 
	sequence with values in a metric space $X$. Let us assume that 
	for every $m\in\n$ there exists
	$$A_{m,\infty}:=\lim_{n\to +\infty}A_{m,n},$$
	and that there exists
	$$A_{\infty,\infty}:=\lim_{m\to +\infty}A_{m,\infty}.$$
	
	Then we have the following conclusions.
	\begin{enumerate}
		\renewcommand{\labelenumi}{(\arabic{enumi})} 
		\item \emph{(Standard conclusion)} There exists a sequence
		$m_{k}\to +\infty$ of nonnegative integers such that
		\begin{equation}
			\lim_{k\to +\infty}A_{m_{k},k}=A_{\infty,\infty}.			
			\label{th:amn-limit}
		\end{equation}
	
		\item \emph{(Refined conclusion)} For every sequence of real
		numbers $r_k\to +\infty$ there exists a sequence $m_{k}\to
		+\infty$ of nonnegative integers such that $m_{k}\leq r_{k}$ 
		for every $k$ large enough, and such that 
		(\ref{th:amn-limit}) holds true.
		\qed
	\end{enumerate}
\end{lemma}

In the next result we estimate from below the norm of a discrete
derivative.  In the continuous setting, when we know the values of
some function $f\in H^{1}_{0}((0,1))$ in some given points, then we
can estimate $\|f_{x}\|_{2}$ from below.  The conclusions of the
following lemma are the discrete counterpart of such estimates.

\begin{lemma}[Discrete derivative estimates]\label{lemma:triang}
	Let $n$ be a positive integer. Let $f\in\PCn$ be a function which 
	is equal to 0 in the last subinterval $(1-1/n,1)$. Let us 
	consider the discrete derivative $D^{-1/n}f(x)$, defined after 
	setting $f(x)=0$ in $(-1/n,0)$.
	
	Then the following estimates hold true.
	\begin{enumerate}
		\renewcommand{\labelenumi}{(\arabic{enumi})}
		\item  We have that
		\begin{equation}
			\left\|D^{-1/n}f(x)\right\|_{2}\geq 2\|f(x)\|_{\infty}.
			\label{th:triang}
		\end{equation}
	
		\item  Let $k$ be a positive integer, and let 
		$0<d_{1}<\ldots<d_{k}<1$. Then we have that
		\begin{equation}
			\left\|D^{-1/n}f(x)\right\|_{2}^{2}\geq\sum_{h=0}^{k}
			\frac{1}{d_{h+1}-d_{h}+1/n}|f(d_{h+1})-f(d_{h})|^{2},
			\label{th:triang-refined}
		\end{equation}
		with the agreement that $d_{0}=0$, $d_{k+1}=1$, and
		$f(d_{0})=f(d_{k+1})=0$.
	\end{enumerate}
\end{lemma}

\paragraph{\textmd{\emph{Proof}}}

Let $f_{i}$ (with $i=1,\ldots,n$) denote the value of $f$ in the
interval $((i-1)/n,i/n)$.  Let us set $f_{0}:=0$ (this choice is
consistent with our extension of $f(x)$ in $(-1/n,0)$), and let us
recall that $f_{n}=0$ due to our assumption on $f$.  Let
$j\in\{1,\ldots,n\}$ be the index (or one of the indices) such that
$\|f(x)\|_{\infty}=|f_{j}|$.

Then we have that
$$\left\|D^{-1/n}f(x)\right\|_{2}\ \geq\ \left\|D^{-1/n}f(x)\right\|_{1}
\ =\ \sum_{i=1}^{j}|f_{i}-f_{i-1}|+ \sum_{i=j+1}^{n}|f_{i}-f_{i-1}|$$
$$\geq\ \left|\sum_{i=1}^{j}(f_{i}-f_{i-1})\right|+
\left|\sum_{i=j+1}^{n}(f_{i}-f_{i-1})\right|\ =\ 
|f_{j}-f_{0}|+|f_{n}-f_{j}|\ =\ 2|f_{j}|,$$
which proves (\ref{th:triang}).

Let us consider now the second statement.  Let us set $i_{0}:=0$,
$i_{k+1}:=n$, and $i_{h}:=\lfloor nd_{h}\rfloor+1$ for every
$h=1,\ldots,k$.  With this notation we have that $f(d_{h})=f_{i_{h}}$
for every $h=0,1,\ldots,k+1$.  Moreover it is easy to see that
$$\frac{i_{h+1}}{n}-\frac{i_{h}}{n}= \frac{\lfloor
nd_{h+1}\rfloor+1}{n}-\frac{\lfloor nd_{h}\rfloor+1}{n}
\leq d_{h+1}-d_{h}+\frac{1}{n}.$$

Thus from H\"{o}lder's inequality it follows that
\begin{eqnarray}
	\left	\|D^{-1/n}f(x)\right\|_{\LL((0,1))}^{2} & = & \sum_{h=0}^{k}
	\left\|D^{-1/n}f(x)\right\|_{L^{2}((i_{h}/n,i_{h+1}/n))}^{2}
	\nonumber  \\
	 & \geq & \sum_{h=0}^{k}
	 \frac{1}{(i_{h+1}/n)-(i_{h}/n)}
	\|D^{-1/n}f(x)\|_{L^{1}((i_{h}/n,i_{h+1}/n))}^{2}
	\nonumber  \\
	 & \geq & \sum_{h=0}^{k}
	 \frac{1}{d_{h+1}-d_{h}+1/n}
	\|D^{-1/n}f(x)\|_{L^{1}((i_{h}/n,i_{h+1}/n))}^{2}.
	\label{est:triang1}
\end{eqnarray}

To be precise, the first inequality requires that all indices $i_{h}$
are distinct, and this is true only when $n$ is large enough.  On the
other hand, the final result is true in any case, because it is enough
to ignore the intervals of length zero in the first sum.  

Finally we have that
$$\left\|D^{-1/n}f(x)\right\|_{L^{1}((i_{h}/n,i_{h+1}/n))}=
\sum_{i=i_{h}+1}^{i_{h+1}}|f_{i}-f_{i-1}|\geq
\left|\sum_{i=i_{h}+1}^{i_{h+1}}(f_{i}-f_{i-1})
\right|=$$
$$=|f_{i_{h+1}}-f_{i_{h}}|=
|f(d_{h+1})-f(d_{h})|.$$

From the last estimate and (\ref{est:triang1}) we obtain 
(\ref{th:triang-refined}).
\qed
\medskip

In the next statement $n$ is fixed, and we present several estimates
relating $k$-energies, the norm of their gradient (the slope), jump
heights, and subcritical incremental quotients.  These estimates are
the technical core of our analysis.

\begin{lemma}[Fundamental estimates]\label{lemma:fund-est}
	Let $k$ and $n$ be positive integers, let $D\subseteq(0,1)$ be a
	finite set with $|D|=k$, and let $v\in\PSDn$.  Let $D_{n}$ be
	defined according to (\ref{defn:Dn}), and let us assume that
	$n$ is big enough so that $(0,1)\setminus D_{n}$ has $(k+1)$
	connected components.
	
	Let $\Gkn(v)$ be the functional defined in (\ref{defn:gkn}), let
	$\nabla\Gkn(v)$ be its gradient, let $J_{d,n}(v)$ be the discrete jump
	heights of $v$ defined in (\ref{defn:Jdn}), and let $SQ_{n}(v)$ be
	the subcritical incremental quotient of $v$ defined in
	(\ref{defn:SQn}).
	
	Then we have that
	\begin{eqnarray}
		 & \displaystyle{\Gkn(v)\geq k\log
		 \left(\min_{d\in D}|J_{d,n}(v)|\right),} & 
		\label{est:f-vs-jh-min}  \\
		\noalign{\vspace{1ex}}
		 &
		 \displaystyle{\Gkn(v)\leq\frac{n}{2}\log
		 \left(1+\left[SQ_{n}(v)\right]^{2}\right)+
		 \frac{1}{2}\sum_{d\in D}\log\left(
		 \frac{1}{n^{2}}+\left[J_{d,n}(v)\right]^{2}\right),} &
		\label{est:f-above}  \\
		\noalign{\vspace{1ex}}
		 & \displaystyle{\left\|\nabla\Gkn(v)\right\|_{2}\geq
			n\left|SQ_{n}(v)\right|,} & 
		\label{est:slope-vs-sq}  \\
		\noalign{\vspace{1ex}}
		 & \displaystyle{\left\|\nabla\Gkn(v)\right\|_{2}\geq
			\left(\min_{d\in D}|J_{d,n}(v)|\right)^{-1}.} & 
		\label{est:s-vs-jh-min2}
	\end{eqnarray}
	
	Finally, if $D=\{d_{1},\ldots,d_{k}\}$ with
	$0<d_{1}<\ldots<d_{k}<1$, then we have that
	\begin{equation}
		\left\|\nabla\Gkn(v)\right\|_{2}^{2}\geq\sum_{i=0}^{k}
		\frac{1}{d_{i+1}-d_{i}+n^{-1}}\left(
		\frac{J_{d_{i+1},n}(v)}{n^{-2}+[J_{d_{i+1},n}(v)]^{2}}-
		\frac{J_{d_{i},n}(v)}{n^{-2}+[J_{d_{i},n}(v)]^{2}}\right)^{2},
		\label{est:s-vs-jhs}
	\end{equation}
	with the agreement that $d_{0}=0$, $d_{k+1}=1$, and 
	$J_{d_{0},n}(v)=J_{d_{k+1},n}(v)=0$.
\end{lemma}

\paragraph{\textmd{\emph{Proof of estimates on the functional}}}

From the definition of jump heights we have that
\begin{eqnarray}
	\frac{n}{2}\int_{D_{n}}
	\log\left(1+|D^{1/n}v(x)|^{2}\right)\,dx-k\log n & = &
	\frac{1}{2}\sum_{d\in D}\log\left(1+n^{2}
	\left[J_{d,n}(v)\right]^{2}\right)-
	\frac{k}{2}\log n^{2}
	\nonumber  \\
	 & = & \frac{1}{2}\sum_{d\in D}\log\left(\frac{1}{n^{2}}+
	 \left[J_{d,n}(v)\right]^{2}\right).
	\label{eq:gkn-dn}
\end{eqnarray}

In order to prove (\ref{est:f-vs-jh-min}), we estimate the right-hand side
of (\ref{defn:gkn}) from below by considering only the integration
over $D_{n}$.  From (\ref{eq:gkn-dn}) we deduce that
$$\Gkn(v) \geq \frac{n}{2}\int_{D_{n}}
\log\left(1+|D^{1/n}v(x)|^{2}\right)\,dx-k\log n =
\frac{1}{2}\sum_{d\in D}\log\left(\frac{1}{n^{2}}+
\left[J_{d,n}(v)\right]^{2}\right)$$
$$\geq\sum_{d\in D}\log|J_{d,n}(v)|\geq 
k\log\left(\min_{d\in D}|J_{d,n}(v)|\right),$$
which proves (\ref{est:f-vs-jh-min}).

On the other hand, from (\ref{defn:SQn}) we have that
\begin{equation}
	\frac{n}{2}\int_{[0,1]\setminus D_{n}}
	\log\left(1+|D^{1/n}v(x)|^{2}\right)\,dx\leq
	\frac{n}{2}\log\left(1+[SQ_{n}(v)]^{2}\right).
	\label{eq:gkn-dnc}
\end{equation}

Summing (\ref{eq:gkn-dnc}) and (\ref{eq:gkn-dn}) we obtain 
(\ref{est:f-above}).

\paragraph{\textmd{\emph{Proof of estimates on the slope}}}

Let us consider the function
\begin{equation}
	f(x):=n\,\frac{D^{1/n}v(x)}{1+\left[D^{1/n}v(x)\right]^{2}}.
	\label{defn:f}
\end{equation}

It turns out that $f(x)=0$ in the last subinterval $(1-1/n,1)$.
Moreover $\nabla\Gkn(v)$, whose expression has been computed in
(\ref{defn:grad-gkn}), coincides (up to the sign) with the discrete
derivative $D^{-1/n}f(x)$, computed after setting $f(x)=0$ in the
interval $(-1/n,0)$.

Therefore, applying (\ref{th:triang}) to the function $f(x)$ defined in
(\ref{defn:f}), we obtain that
\begin{equation}
	\left\|\nabla\Gkn(v)\right\|_{2}=
	\left\|D^{-1/n}f(x)\right\|_{2}\geq 2\|f(x)\|_{\infty}.
	\label{est:slope-triang}
\end{equation}

Let us estimate the right-hand side of (\ref{est:slope-triang}) from
below by restricting the $L^{\infty}$-norm to $[0,1]\setminus D_{n}$.
Since the function $\sigma\to\sigma(1+\sigma^{2})^{-1}$ is increasing
in $[-1,1]$, and since $0\leq SQ_{n}(v)\leq 1$, we obtain that
$$\|f(x)\|_{\infty}\geq\max\left\{
\frac{n\,|D^{1/n}v(x)|}{1+\left[D^{1/n}v(x)\right]^{2}}:
x\in[0,1]\setminus D_{n}\right\}=
\frac{n\,SQ_{n}(v)}{1+\left[SQ_{n}(v)\right]^{2}}
\geq\frac{n}{2}\,SQ_{n}(v),$$
which proves (\ref{est:slope-vs-sq}). 

Now let us estimate the right-hand side of (\ref{est:slope-triang}) from
below by restricting the $L^{\infty}$-norm to $D_{n}$. We obtain that
$$\|f(x)\|_{\infty}\geq\max\left\{
\frac{n\,|D^{1/n}v(x)|}{1+\left[D^{1/n}v(x)\right]^{2}}:
x\in D_{n}\right\}=\max_{d\in D}
\frac{n^{2}\,|J_{d,n}(v)|}{1+n^{2}\left[J_{d,n}(v)\right]^{2}}.$$

Since $n|J_{d,n}(v)|\geq 1$ for
every $d\in D$, and since the function
$\sigma\to\sigma(1+\sigma^{2})^{-1}$ is decreasing for
$\sigma\geq 1$, we have that the maximum is achieved when the argument
is minimum, hence
$$\left\|\nabla\Gkn(v)\right\|_{2}\geq 2\|f(x)\|_{\infty}\geq
2\max_{d\in D}
\frac{n^{2}\,|J_{d,n}(v)|}{1+n^{2}\left[J_{d,n}(v)\right]^{2}}=
\frac{2n^{2}\,\displaystyle{\min_{d\in D}|J_{d,n}(v)|}}{1+n^{2}
\Bigl[\displaystyle{\min_{d\in D}|J_{d,n}(v)|}\Bigl]^{2}}.$$

Recalling once more that the minimum is greater than or equal to
$1/n$, estimate (\ref{est:s-vs-jh-min2}) follows.

A more refined estimate, keeping into account all jump heights,
follows from (\ref{th:triang-refined}) applied to the function $f(x)$
defined in (\ref{defn:f}).  Since $D^{1/n}v(d)=nJ_{d,n}(v)$ for every
$d\in D$, we obtain exactly (\ref{est:s-vs-jhs}).\qed 
\medskip

In the next lemma we consider the difference between two piecewise
constant functions $v$ and $w$.  The main idea is the following.  If
the number of discrete jump points is finite, and their location is
fixed, then the $L^{2}$-norm of $v-w$ estimates both the
$L^{\infty}$-norm of $v-w$, and the difference between jump heights.

\begin{lemma}[Uniform and jump-height estimates]\label{lemma:jh-est}
	Let $D\subseteq(0,1)$ be a finite set, let $K_{0}$ be the minimum
	of the lengths of the $(|D|+1)$ intervals into which $(0,1)$ is
	divided by $D$, and let $n\geq 3K_{0}^{-1}$ be a positive integer.
	
	Let $D'\subseteq D$ and $D''\subseteq D$ be two subsets, and let
	$v\in PS_{D',n}$ and $w\in PS_{D'',n}$ be two piecewise constant
	functions with discrete jump heights $J_{d,n}(v)$ and
	$J_{d,n}(w)$, respectively (we can think both jump heights as
	defined for every $d\in D$).  
	
	Then we have that
	\begin{equation}
		\min\left\{\strut
		K_{0},\left\|v-w\right\|_{\infty}\right\}\leq
		3\|v-w\|_{2}^{2/3},
		\label{th:est-i-2}
	\end{equation}
	\begin{equation}
		\min\left\{\strut
		K_{0},\left|J_{d,n}(v)-J_{d,n}(w)\right|\right\}\leq
		6\|v-w\|_{2}^{2/3}
		\quad\quad
		\forall d\in D.
		\label{th:est-d-d}
	\end{equation}
\end{lemma}

\paragraph{\textmd{\emph{Proof}}}

Let $f(x):=|v(x)-w(x)|$, and let $f_{i}$ (with $i=1,\ldots,n$) denote
the value of $f$ in the $i$-th subinterval.  Let $j\in\{1,\ldots,n\}$
be the index (or one of the indices) such that
$\|f(x)\|_{\infty}=f_{j}$.

For every $x\in[0,1]\setminus D_{n}$ we have that $|D^{1/n}v(x)|\leq
1$ and $|D^{1/n}w(x)|\leq 1$, hence
\begin{equation}
	|D^{1/n}f(x)|\leq 2
	\quad\quad
	\forall x\in[0,1]\setminus D_{n}.
	\label{est:f-2lip}
\end{equation}

Let us set $H:=\min\{K_{0},f_{j}\}$, and let us consider the two 
intervals
$$I:=\left(\frac{j}{n}-\frac{H}{3},\frac{j}{n}\right),
\hspace{3em}
I':=\left(\frac{j-1}{n},\frac{j-1}{n}+\frac{H}{3}\right).$$

Since $n\geq 3K_{0}^{-1}$, estimate (\ref{est:f-2lip}) implies that 
in at least one of these intervals the difference between the values of $f$ in 
any two neighboring subintervals is always less than or equal to 
$2/n$. Let us assume, without loss of generality, that this
happens in $I$ (the other case is specular). Then we have that
$$f(x)\ \geq\ f_{j}-2\left(\frac{j}{n}-x\right)\ \geq\ 
H-2\,\frac{H}{3}\ \geq\ 
\frac{H}{3} \quad\quad
\forall x\in I,$$
so that
$$\|v-w\|^{2}_{L^{2}(I)}=\int_{I}[f(x)]^{2}\,dx\geq\frac{H^{3}}{27},$$
and therefore 
$$\min\left\{K_{0},\|v-w\|_{\infty}\right\}=H\leq
3\|v-w\|^{2/3}_{L^{2}(I)}\leq
3\|v-w\|^{2/3}_{\LL((0,1))}.	$$
	
This proves (\ref{th:est-i-2}). Now we have that
$$|J_{d,n}(v)-J_{d,n}(w)|=|J_{d,n}(v-w)|\leq 2\|v-w\|_{\infty}
\quad\quad
\forall d\in D,$$
hence 
$$\min\left\{K_{0},\left|J_{d,n}(v)-J_{d,n}(w)\right|\strut\right\}\leq
2 \min\left\{K_{0},\|v-w\|_{\infty}\right\}\leq
6\|v-w\|_{2}^{2/3}$$
for every $d\in D$, which is exactly (\ref{th:est-d-d}).\qed
\medskip

In the last lemma we consider a sequence $v_{n}\to v$.  We point out
that this sequence is allowed to loose jump points in the limit.

\begin{lemma}[Jump convergence and BV estimate]\label{lemma:unif-est}
	Let $D'\subseteq D\subseteq (0,1)$ be two finite sets.  Let $v\in
	PS_{D'}$, and let $\{v_{n}\}$ be a sequence such that $v_{n}\to 
	v$ in $\LL$.
	Let us assume that for every $n\geq 1$ we have that $v_{n}\in 
	PS_{D''(n),n}$ for some finite set $D''(n)\subseteq D$.
	
	Then we have the following conclusions.
	\begin{enumerate}
		\renewcommand{\labelenumi}{(\arabic{enumi})}
		\item\label{stat:jc} \emph{(Jump convergence)} Let $J_{d}(v)$
		be the jump heights of $v$, and let $J_{d,n}(v_{n})$ be the
		discrete jump heights of $v_{n}$ (we can think both jump
		heights as defined for every $d\in D$, with the agreement that
		$J_{d}(v)=0$ when $d\in D\setminus D'$).  Then we have
		that
		\begin{equation}
			\lim_{n\to +\infty}J_{d,n}(v_{n})=J_{d}(v)
			\quad\quad
			\forall d\in D.
			\label{th:jump-conv}
		\end{equation}
		As a consequence, $D'\subseteq D''(n)$ for every $n$ large 
		enough.
		
		\item\label{stat:unif-est} \emph{(Uniform $BV$ estimates)} We have that 
		\begin{equation}
			\sup_{n\geq 1} \left\{\left\|D^{1/n}v_{n}\right\|_{1}+
			\|v_{n}\|_{\infty}\right\}<+\infty.
			\label{th:tv-li-unif}
		\end{equation}
	\end{enumerate}
\end{lemma}

\paragraph{\textmd{\emph{Proof of statement (\ref{stat:jc})}}}

Let $w_{n}\in\PCn$ be the piecewise constant approximation of $v$
defined by 
$$w_{n}(x):=v\left(\frac{\lfloor nx\rfloor}{n}\right)
\quad\quad
\forall x\in [0,1],\ \forall n\geq 1.$$

Let $K_{0}$ be as in Lemma~\ref{lemma:jh-est}, and let $n\geq
3K_{0}^{-1}$ as in that lemma.  It is not difficult to see that
$$\left|J_{d}(v)-J_{d,n}(w_{n})\right|\leq 2n^{-1}  \quad\quad
\forall d\in D,$$
hence
\begin{eqnarray*}
	\left|J_{d}(v)-J_{d,n}(v_{n})\right| & \leq &
	\left|J_{d}(v)-J_{d,n}(w_{n})\right|+
	\left|J_{d,n}(w_{n})-J_{d,n}(v_{n})\right| \\
	 & \leq &
	2n^{-1}+\left|J_{d,n}(w_{n})-J_{d,n}(v_{n})\right|.
\end{eqnarray*}

Therefore, applying (\ref{th:est-d-d}) with $v=v_{n}$ and 
$w=w_{n}$, we obtain that
\begin{eqnarray*}
	\min\left\{K_{0},|J_{d}(v)-J_{d,n}(v_{n})|\strut\right\} & \leq & 
	\min\left\{K_{0},|J_{d,n}(w_{n})-J_{d,n}(v_{n})|
	+2n^{-1}\right\}\\
	 & \leq & \min\left\{\strut K_{0},|J_{d,n}(w_{n})-J_{d,n}(v_{n})|
	\right\}+2n^{-1}  \\
	 & \leq & 6\|w_{n}-v_{n}\|_{2}^{2/3}+2n^{-1}  \\
	 & \leq & 6\left(\|w_{n}-v\|_{2}+
	 \|v-v_{n}\|_{2}\strut\right)^{2/3}+2n^{-1}
\end{eqnarray*}
for every $d\in D$.  All the terms in the right-hand side tend to 0 as
$n\to +\infty$.  This proves~(\ref{th:jump-conv}).

\paragraph{\textmd{\emph{Proof of statement (\ref{stat:unif-est})}}}

Let $D_{n}$ be defined by (\ref{defn:Dn}).
For every $d\in D$ we have that
$D^{1/n}v_{n}(x)=nJ_{d,n}(v_{n})$ for every $x$ in the corresponding
subinterval, hence $$\left\|D^{1/n}v_{n}\right\|_{1}=
\int_{D_{n}}|D^{1/n}v_{n}(x)|\,dx+ \int_{[0,1]\setminus
D_{n}}|D^{1/n}v_{n}(x)|\,dx\leq
\sum_{d\in D}|J_{d,n}(v_{n})|+1.$$

Due to (\ref{th:jump-conv}) we have that
$$\lim_{n\to +\infty}\sum_{d\in D}|J_{d,n}(v_{n})|+1=
\sum_{d\in D'}|J_{d}(v)|+1<+\infty,$$
which is enough to prove the equi-boundedness of the total variations
$\|D^{1/n}v_{n}\|_{1}$.

The uniform bound on $\|v_{n}\|_{\infty}$ follows from the uniform
bound on total variations, and from the fact that the average of
$v_{n}$ tends to the average of $v$ (hence averages are
equi-bounded).\qed \medskip

\subsection{Evolution problems as maximal slope curves}

Both the differential equation (\ref{pbm:main-eq}), and the system of
ordinary differential equations introduced in
section~\ref{sec:limit-pbm}, are equivalent to suitable integral
(in)equalities.  This equivalence is the key point in the theory of
maximal slope curves, for which we refer to~\cite{AGS}.

For the sake of simplicity, we want to keep this paper as independent
as possible of the general abstract theory of maximal slope curves.
For this reason, in Proposition~\ref{prop:cmp} we state the two
implications we need throughout this paper, and we provide a
self-contained and almost elementary proof of them.  Before stating
these implications, we need the following definition.

\begin{defn}[Slope of limit functional]\label{defn:slope-gki}
	\begin{em}
		Let $D\subseteq(0,1)$ be a finite set with $|D|=k$.  Let us
		assume that $D=\{d_{1},\ldots,d_{k}\}$, with
		$0<d_{1}<\ldots<d_{k}<1$.  Let $v\in\PCD$, and let $J_{d}(v)$
		be its jump heights.  Let $\Gki$ be the functional defined in
		(\ref{defn:gki}).
		
		The \emph{slope} of $\Gki$ in the point $v$ with respect to
		the $\LL$-metric is the nonnegative real number whose square is
		given by
		\begin{eqnarray*}
			\left\|\nabla\Gki(v)\right\|_{2}^{2} & := &
			\frac{1}{d_{1}}\cdot\frac{1}{[J_{d_{1}}(v)]^{2}}+
			\frac{1}{1-d_{k}}\cdot\frac{1}{[J_{d_{k}}(v)]^{2}}
			\\
			 &  & 
			 +\sum_{i=1}^{k-1}
			\frac{1}{d_{i+1}-d_{i}}\left( \frac{1}{J_{d_{i+1}}(v)}-
			\frac{1}{J_{d_{i}}(v)}\right)^{2}.
		\end{eqnarray*}
	\end{em}
\end{defn}

There are several interpretations of the slope. In this case the 
domain where the functional is finite is isometric to an open subset 
of $\re^{k+1}$, and under this isometry the functional can be 
identified with a function of $(k+1)$ real variables. The slope 
coincides with the norm of the gradient of this function. 

The approximated solution $u_{n}(t)$ is the gradient-flow in $\LL$ of
the functional $\Gkn$ for every $t\geq 0$, while the function $u(t)$
defined in section~\ref{sec:limit-pbm} is the gradient-flow in $\LL$ of
the functional $\Gki$ up to $\tsing$.  In any case, what we need in
this paper (and in particular in the proof of
Proposition~\ref{prop:fino-t1}) are the following two facts.

\begin{prop}[Differential equations vs maximal slope
curves]\label{prop:cmp} 
\hskip 0em plus 2 em 
The approximating problems and the limit problem can be reformulated 
as follows.
\begin{enumerate}
	\renewcommand{\labelenumi}{(\arabic{enumi})}
	\item \emph{(Approximating problems)} Let $n$ be a positive
	integer, and let $u_{0n}\in\PCn$.  Let $u_{n}\in
	C^{1}\left([0,+\infty);\LL\right)$ be the solution of
	(\ref{pbm:main-eq}), (\ref{pbm:main-data}).  Then we have that
	\begin{equation}
		\Gkn(u_n(s))-\Gkn(u_n(t))= \frac 1 2
		\int_s^t \|u_n'(\tau)\|_{2}^2\,d\tau+ 
		\frac 1 2 \int_s^t \left\|\nabla \Gkn
		(u_n(\tau))\right\|_{2}^2\,d\tau
		\label{hp:cmp-gkn}
	\end{equation}
	for every $0\leq s\leq t$ and every $k\in\n$.

	\item \emph{(Limit problem)} Let $k$ be a positive integer,
	let $D\subseteq(0,1)$ be a finite set with $|D|=k$, and let
	$u_{0}\in\PCD$.  Let $T_{0}>0$, and let $v\in
	H^{1}\left((0,T_{0});\LL\right)$ be a function such that
	$v(0)=u_{0}$, $v(t)\in\PCD$ for every $t\in[0,T_{0}]$, and
	\begin{equation}
		\Gki(v(s))-\Gki(v(t))\geq \frac 1 2 \int_s^t
		\|v'(\tau)\|_{2}^2\,d\tau+ \frac 1 2 \int_s^t \left\|\nabla \Gki
		(v(\tau))\right\|_{2}^2\,d\tau			
		\label{hp:cmp-gki}
	\end{equation}

	for every $0< s\leq t< T_{0}$.  Then $v(t)$ coincides 
	in $[0,T_{0}]$ with the function $u(t)$ defined in 
	section~\ref{sec:limit-pbm}.
	
\end{enumerate}
\end{prop}

\paragraph{\textmd{\emph{Proof}}}

From equation (\ref{pbm:main-eq}), which is the same as
(\ref{pbm:main-eq-gkn}), we have that 
$$-\frac{d}{dt}\Gkn(u_{n}(t))=
-\langle\nabla\Gkn(u_{n}(t)),u_{n}'(t)\rangle=
\frac{1}{2}\|u_{n}'(t)\|_{2}^{2}+
\frac{1}{2}\left\|\Gkn(u_{n}(t))\right\|_{2}^{2}.$$

Integrating in $[s,t]$ we obtain (\ref{hp:cmp-gkn}).

Let us consider now inequality (\ref{hp:cmp-gki}).  Let
$D=\{d_{1},\ldots,d_{k}\}$ with $0<d_{1}<\ldots<d_{k}<1$, let
$d_{0}:=0$ and $d_{k+1}:=1$, and let us identify $v(t)$ with the
vector of plateau heights $(a_{0}(t),a_{1}(t),\ldots,a_{k}(t))$, where
$a_{i}(t)$ is the constant value of $v(t)$ for $x\in(d_{i},d_{i+1})$.

The $H^{1}$ regularity of $v(t)$ implies $H^{1}$ regularity of all 
components. Let us compute the time derivative of the function 
$t\to\Gki(v(t))$. Using the chain rule, and rearranging the terms, 
for almost every $t\in(0,T_{0})$ we obtain that
\begin{eqnarray*}
	\lefteqn{\hspace{-3em}-\frac{d}{dt}\Gki(v(t))\ =\
	-\frac{d}{dt}\sum_{i=1}^{k}\log|a_{i}(t)-a_{i-1}(t)|} \\
	\quad\quad &  = &
	\sum_{i=0}^{k}a_{i}'(t) \left(\frac{1}{a_{i+1}(t)-a_{i}(t)}-
	\frac{1}{a_{i}(t)-a_{i-1}(t)}\right)
	\\
	\quad\quad &  = & \sum_{i=0}^{k}\sqrt{d_{i+1}-d_{i}}\,a_{i}'(t)\cdot
	 \frac{1}{\sqrt{d_{i+1}-d_{i}}}
	 \left(\frac{1}{a_{i+1}(t)-a_{i}(t)}-
	\frac{1}{a_{i}(t)-a_{i-1}(t)}\right),
\end{eqnarray*}
with the agreement to neglect the two fractions involving indices 
less than 0 or larger than $k$ (which appear in the terms of the sum 
corresponding to $i=0$ and $i=k$).

Applying the inequality $xy\leq 2^{-1}(x^{2}+y^{2})$ to each term of 
the sum, we find that (for shortness' sake we drop the dependence on $t$ in 
the right-hand side of the first line)
\begin{eqnarray}
	-\frac{d}{dt}\Gki(v(t)) &  \leq & 
	 \frac{1}{2}\sum_{i=0}^{k}(d_{i+1}-d_{i})\left[a_{i}'\right]^{2}
	 +\frac{1}{2}\sum_{i=0}^{k}\frac{1}{d_{i+1}-d_{i}}\left(
	 \frac{1}{a_{i+1}-a_{i}}-\frac{1}{a_{i}-a_{i-1}}\right)^{2}
	\nonumber  \\
	 & = & \frac{1}{2}\|v'(t)\|_{2}^{2}+\frac{1}{2}
	 \left\|\nabla\Gki(v(t))\right\|_{2}^{2}
	\label{eqn:gki-leq}
\end{eqnarray}
for almost every $t\in(0,T_{0})$. On the other hand, from 
(\ref{hp:cmp-gki}) we easily obtain that
\begin{equation}
-\frac{d}{dt}\Gki(v(t))\geq\frac{1}{2}\|v'(t)\|_{2}^{2}+\frac{1}{2}
 \left\|\nabla\Gki(v(t))\right\|_{2}^{2}
\label{eqn:gki-geq}
\end{equation}
for almost every $t\in(0,T_{0})$.
Comparing (\ref{eqn:gki-leq}) and (\ref{eqn:gki-geq}) we deduce that 
there is equality for almost every $t\in(0,T_{0})$. But in the 
inequality used to deduce (\ref{eqn:gki-leq}) we have equality if and only if
$$\sqrt{d_{i+1}-d_{i}}\ a_{i}'(t)=
\frac{1}{\sqrt{d_{i+1}-d_{i}}}\left(
\frac{1}{a_{i+1}(t)-a_{i}(t)}-\frac{1}{a_{i}(t)-a_{i-1}(t)}\right)$$
for every $i=0,\ldots,k$ and for almost every $t\in(0,T_{0})$, which
is equivalent to the system of ordinary differential equations
introduced in section~\ref{sec:limit-pbm}.  Since the right-hand side
is continuous, we can conclude that actually $a_{i}(t)$ is of class
$C^{1}$, and we have equality for every $t$ in the closed interval
$[0,T_{0}]$.\qed \medskip

The formulation in terms of differential inequalities is very stable
when passing to the limit.  In the proof of
Proposition~\ref{prop:fino-t1} we obtain (\ref{hp:cmp-gki}) by passing
to the limit in (\ref{hp:cmp-gkn}) as $n\to +\infty$.  The following
result is fundamental in that stage.

\begin{prop}[Bounded slope sequences]\label{prop:cmp-limit}
	\hskip 0em plus 1 em
	Let $k$ be a nonnegative integer, and let $D\subseteq (0,1)$ be a
	finite set with $|D|=k$.  Let $v\in\LL$, and let $\{v_{n}\}$ be a
	sequence such that $v_{n}\in\PSDn$ for every $n\geq 1$, and such 
	that $v_{n}\to v$ in $\LL$.
	
	Let $J_{d,n}(v_{n})$ denote the discrete jump heights of $v_{n}$,
	and let us suppose that there exist $c_{0}>0$ and $n_{0}\geq 1$
	such that
	\begin{equation}
		|J_{d,n}(v_{n})|\geq c_{0}
		\quad\quad
		\forall d\in D,\ \forall n\geq n_{0}.		
		\label{hp:jh-big}
	\end{equation}
	
	Then we have the following conclusions.
	\begin{enumerate}
		\renewcommand{\labelenumi}{(\arabic{enumi})}
		\item  \emph{(Gamma-liminf inequality for slopes)} We 
		have that
		\begin{equation}
			\liminf_{n\to +\infty}\left\|\nabla\Gkn(v_{n})\right\|_{2}\geq
			\left\|\nabla\Gki(v)\right\|_{2},
			\label{th:slope-liminf}
		\end{equation}
		where the right-hand side is intended to be $+\infty$ if 
		$v\not\in\PCD$.
	
		\item  \emph{(Bounded slope sequences are recovery 
		sequences)} If in addition we assume that
		\begin{equation}
			\sup_{n\geq 1}\left\|\nabla\Gkn(v_{n})\right\|_{2}<+\infty,
			\label{hp:bounded-slope}
		\end{equation}
		then we have that $v\in\PCD$, and moreover
		\begin{eqnarray}
			 & |J_{d}(v)|\geq c_{0}
			 \quad\quad
			 \forall d\in D, & 
			\label{th:jh-big}  \\
			\noalign{\vspace{1ex}}
			 & \displaystyle{\lim_{n\to +\infty}\Gkn(v_{n})=\Gki(v).} & 
			\label{th:bss-cont}
		\end{eqnarray}
	\end{enumerate}
\end{prop}

\paragraph{\textmd{\emph{Proof}}}

We prove the two statements in reverse order.

\subparagraph{\textmd{\emph{Statement (2)}}}

Let $c_{1}$ denote the supremum in (\ref{hp:bounded-slope}).  From
(\ref{est:slope-vs-sq}) we have that
\begin{equation}
	c_{1}\geq\left\|\nabla\Gkn(v_{n})\right\|_{2}\geq
	n\,SQ_{n}(v_{n})
	\quad\quad
	\forall n\geq 1.
	\label{est:siq}
\end{equation}

Let $\delta>0$, and let $D_{\delta}$ denote the neighborhood of $D$ 
with width $\delta$. Then (\ref{est:siq}) implies that
$$D^{1/n}v_{n}(x)\to 0
\quad\quad
\mbox{uniformly in }[0,1]\setminus D_{\delta},$$
which in turn implies that $v(x)$ is constant in each connected
component of $[0,1]\setminus D_{\delta}$.

Since $\delta$ is arbitrary, this proves that $v\in PC_{D'}$ for some
$D'\subseteq D$.  On the other hand, assumption (\ref{hp:jh-big}) and
the convergence of jump heights (\ref{th:jump-conv}) imply
(\ref{th:jh-big}), which proves also that actually $D'=D$.

In order to prove (\ref{th:bss-cont}), we split the integral in 
(\ref{defn:gkn}) into an integral in $D_{n}$, and an integral in 
$[0,1]\setminus D_{n}$. For the second one we apply (\ref{est:siq}) 
and we deduce that
$$0\leq\frac{n}{2}\int_{[0,1]\setminus D_{n}}
\log\left(1+|D^{1/n}v_{n}(x)|^{2}\right)\,dx\leq
\frac{n}{2}\log\left(1+|SQ_{n}(v_{n})|^{2}\right)\leq
\frac{n}{2}\log\left(1+\frac{c_{1}^{2}}{n^{2}}\right),$$
which proves that the integral in $[0,1]\setminus D_{n}$ tends to 0. 

For the integral in $D_{n}$ we apply (\ref{eq:gkn-dn}) to the function
$v_{n}$, and we deduce that 
$$\frac{n}{2}\int_{D_{n}}
\log\left(1+|D^{1/n}v_{n}(x)|^{2}\right)\,dx-k\log n =
\frac{1}{2}\sum_{d\in D}\log\left(\frac{1}{n^{2}}+
\left[J_{d,n}(v_{n})\right]^{2}\right).$$

From the jump convergence (\ref{th:jump-conv}) it follows that this
expression tends to
$$\sum_{d\in D}\log|J_{d}(v)|=\Gki(v),$$
which completes the proof of (\ref{th:bss-cont}).

\subparagraph{\textmd{\emph{Statement (1)}}}

Let us take any subsequence (not relabeled) which realizes the
$\liminf$ in (\ref{th:slope-liminf}).  We can assume that
(\ref{hp:bounded-slope}) holds true on this subsequence (otherwise the
conclusion is trivial).  As we have seen in the proof of
statement~(2), this implies in particular that $v\in\PCD$.

Now let us apply estimate (\ref{est:s-vs-jhs}) to the function
$v_{n}$.  We obtain that
\begin{eqnarray*}
	\left\|\nabla\Gkn(v_{n})\right\|_{2}^{2} & \geq &
	\frac{1}{d_{1}+n^{-1}}\left(
	\frac{J_{d_{1},n}(v_{n})}{n^{-2}+[J_{d_{1},n}(v_{n})]^{2}}\right)^{2}  \\
	 &  & +\frac{1}{1-d_{k}+n^{-1}}\left(
	\frac{J_{d_{k},n}(v_{n})}{n^{-2}+[J_{d_{k},n}(v_{n})]^{2}}\right)^{2}  \\
	 &  & +\sum_{i=1}^{k-1}
	\frac{1}{d_{i+1}-d_{i}+n^{-1}}\left(
	\frac{J_{d_{i+1},n}(v_{n})}{n^{-2}+[J_{d_{i+1},n}(v_{n})]^{2}}-
	\frac{J_{d_{i},n}(v_{n})}{n^{-2}+[J_{d_{i},n}(v_{n})]^{2}}\right)^{2}.
\end{eqnarray*}

Let us finally pass to the limit as $n\to +\infty$.  Thanks to the
jump convergence (\ref{th:jump-conv}), the right-hand side tends to
$\|\nabla\Gki(v)\|_{2}^{2}$ (assumption (\ref{hp:jh-big}) guarantees
that all fractions have a finite limit), as defined in
Definition~\ref{defn:slope-gki}.  \qed

\subsection{H\"{o}lder continuity of approximating and limit problems}
\label{sec:proof-u}

\subsubsection*{Proof of Proposition~\ref{prop:u}}

The existence of a local solution to the system of ordinary 
differential equations is trivial. From now on we identify the vector 
$(a_{0}(t),\ldots,a_{k}(t))$ with the function $u(t)$, and we define 
the points $d_{i}$ ($i=0,1,\ldots,k,k+1$) as in 
section~\ref{sec:limit-pbm}.

\paragraph{\textmd{\emph{$L^{\infty}$ estimate}}}

Let $\psi\in C^{1}(\re)$ be an even convex function, and let us 
set
$$\Psi(t):=\int_{0}^{1}\psi(u(t))\,dx=
\sum_{i=0}^{k}(d_{i+1}-d_{i})\psi(a_{i}(t)).$$

Since $\psi'$ is nondecreasing, with some computations it turns out
that
\begin{equation}
	\Psi'(t)=-\sum_{i=1}^{k}
	\frac{\psi'(a_{i}(t))-\psi'(a_{i-1}(t))}{a_{i}(t)-a_{i-1}(t)}
	\leq 0.
	\label{est:deriv-psi}
\end{equation}

Let us assume now in addition that $\psi(\sigma)=0$ if and only if 
$|\sigma|\leq \|u_{0}\|_{\infty}$. Then we have that $\Psi(0)=0$, 
$\Psi(t)\geq 0$ as soon as it is defined, and $\Psi$ is nonincreasing 
because of (\ref{est:deriv-psi}). It follows that $\Psi(t)=0$ as soon 
as it is defined, which proves (\ref{th:infty-u}). 

\paragraph{\textmd{\emph{Finite time break-down}}}

Let us consider the function
$$S(t):=\sum_{d\in D}|J_{d}(u(t))|=
\sum_{i=1}^{k}|a_{i}(t)-a_{i-1}(t)|.$$

The sign of all jump heights is constant as soon as the solution is 
defined. This implies that $S(t)$ is smooth. Computing the time
derivative, and rearranging the terms, we obtain that (for 
shortness's sake we drop the dependence on $t$ in the second line)
\begin{eqnarray*}
	S'(t) & = & -\frac{1}{d_{1}}\,\frac{1}{|a_{1}(t)-a_{0}(t)|}-
	\frac{1}{1-d_{k}}\,\frac{1}{|a_{k}(t)-a_{k-1}(t)|}  \\
	\noalign{\vspace{1ex}}
	&  & -\sum_{i=1}^{k-1}
	\frac{1}{d_{i+1}-d_{i}}\left(
	\frac{1}{|a_{i+1}-a_{i}|}+
	\frac{1}{|a_{i}-a_{i-1}|}\right)
	\left(1-\frac{a_{i}-a_{i-1}}{|a_{i}-a_{i-1}|}\cdot
	\frac{a_{i+1}-a_{i}}{|a_{i+1}-a_{i}|}\right)\\
	 \noalign{\vspace{1ex}}
	 & \leq & -\frac{1}{d_{1}}\,\frac{1}{|a_{1}(t)-a_{0}(t)|}-
	\frac{1}{1-d_{k}}\,\frac{1}{|a_{k}(t)-a_{k-1}(t)|}   \\ 
	\noalign{\vspace{1ex}}
	 & \leq & -\frac{1}{\|u_{0}\|_{\infty}}.
\end{eqnarray*}

Since $S(t)$ is clearly nonnegative, this implies that the solution 
cannot be global, and also provides an estimate on the life span.

\paragraph{\textmd{\emph{H\"{o}lder continuity up to collision}}}

Let us compute the time derivative of the function $t\to\Gki(u(t))$. 
Rearranging the terms we obtain that
\begin{eqnarray*}
	\lefteqn{\hspace{-2em}-\frac{d}{dt}\Gki(u(t))\ =\
	\sum_{i=0}^{k}a_{i}'(t) \left(\frac{1}{a_{i+1}(t)-a_{i}(t)}-
	\frac{1}{a_{i}(t)-a_{i-1}(t)}\right)}
	\\
	 & \quad = & \sum_{i=0}^{k}\sqrt{d_{i+1}-d_{i}}\,a_{i}'(t)\cdot
	 \frac{1}{\sqrt{d_{i+1}-d_{i}}}
	 \left(\frac{1}{a_{i+1}(t)-a_{i}(t)}-
	\frac{1}{a_{i}(t)-a_{i-1}(t)}\right),
\end{eqnarray*}
with the usual agreement to neglect the two fractions involving terms 
$a_{i}(t)$ with indices less than 0 or larger than $k$.

The two factors in each term of the sum are equal due to the system 
of ordinary differential equations. Therefore, the sum can be 
rewritten both in the form
$$\sum_{i=0}^{k}(d_{i+1}-d_{i})[a_{i}'(t)]^{2}=\|u'(t)\|^{2}_{2},$$
and in the form
$$\sum_{i=0}^{k}\frac{1}{d_{i+1}-d_{i}}
\left(\frac{1}{a_{i+1}(t)-a_{i}(t)}-
\frac{1}{a_{i}(t)-a_{i-1}(t)}\right)^{2}=
\left\|\nabla\Gki(u(t))\right\|^{2}_{2}.$$

As a consequence, we have also that
\begin{equation}
	-\frac{d}{dt}\Gki(u(t))=\left\|\nabla\Gki(u(t))\right\|^{2/3}_{2}
	\|u'(t)\|^{4/3}_{2}.
	\label{est:deriv-gki}
\end{equation}

Now let us consider the function
\begin{equation}
	H(t):=3k \exp\left(\frac{2}{3k}\,\Gki(u(t))\right).
	\label{defn:h}
\end{equation}

We claim that the descent of $H(t)$ estimates the 1/4-H\"{o}lder
constant of $u(t)$.  To this end, we begin by computing the time
derivative of $H(t)$.  From (\ref{est:deriv-gki}) we
have that
\begin{equation}
	-H'(t)=2\exp\left(\frac{2}{3k}\,\Gki(u(t))\right)
	\left\|\nabla\Gki(u(t))\right\|_{2}^{2/3}
	\|u'(t)\|_{2}^{4/3}.
	\label{est:deriv-h}
\end{equation}

Let us estimate the first two terms in the right-hand side. For the 
first one we have that
\begin{equation}
	\exp\left(\frac{2}{3k}\,\Gki(u(t))\right)=
	\left[\prod_{i=1}^{k}|a_{i}(t)-a_{i-1}(t)|\right]^{2/(3k)}\geq
	\min_{i=1,\ldots,k}|a_{i}(t)-a_{i-1}(t)|^{2/3}.
	\label{est:exp}
\end{equation}

Let $j$ be the index (or one of the indices) which realizes the 
minimum. Then from Cauchy-Schwarz inequality we have that
\begin{eqnarray}
	\left\|\nabla\Gki(u(t))\right\|^{2}_{2} & \geq & 
	\sum_{i=0}^{j-1}(d_{i+1}-d_{i})\cdot
	\sum_{i=0}^{j-1}\frac{1}{d_{i+1}-d_{i}}
	\left(\frac{1}{a_{i+1}(t)-a_{i}(t)}-
	\frac{1}{a_{i}(t)-a_{i-1}(t)}\right)^{2}
	\nonumber  \\
	 & \geq & \left[\sum_{i=0}^{j-1}
	\left(\frac{1}{a_{i+1}(t)-a_{i}(t)}-
	\frac{1}{a_{i}(t)-a_{i-1}(t)}\right)\right]^{2}
	\nonumber  \\
	 & = & \left(\frac{1}{a_{j}(t)-a_{j-1}(t)}\right)^{2}.
	\label{est:slope}
\end{eqnarray}

Plugging (\ref{est:exp}) and (\ref{est:slope}) into (\ref{est:deriv-h}) we 
obtain that
$$-H'(t)\geq 2\|u'(t)\|_{2}^{4/3}
\geq \|u'(t)\|_{2}^{4/3}
\quad\quad
\forall t\in[0,\tsing).$$

Now we integrate in $[s,t]$, and we exploit that $H(t)$ is nonnegative
and nonincreasing (which follows from (\ref{est:deriv-h})).  We deduce
that 
$$\int_{s}^{t}\|u'(\tau)\|_{2}^{4/3}d\tau\leq H(s)-H(t)\leq H(0)=
3k\exp\left(\frac{2}{3k}\Gki(u_{0})\right)$$
for every $0\leq s\leq t<\tsing$.
Finally, from H\"{o}lder's inequality we obtain that
$$\|u(t)-u(s)\|_{2} \leq
\int_{s}^{t}\|u'(\tau)\|_{2}\,d\tau \leq \left(
\int_{s}^{t}\|u'(\tau)\|_{2}^{4/3}d\tau\right)^{3/4}
|t-s|^{1/4}$$
for every $0\leq s\leq t<\tsing$.  Combining the last two
estimates we obtain (\ref{th:h-1/4-u}).\qed

\subsubsection*{Proof of Proposition~\ref{prop:hc}}

In analogy with (\ref{defn:h}), let us consider the function
$$H_{n}(t) := 3k \exp\left(\frac{2}{3k}\,\Gkn(u_{n}(t))\right).$$

Exploiting equation (\ref{pbm:main-eq-gkn}) (which is the same as
(\ref{pbm:main-eq})), in analogy with (\ref{est:deriv-h}) we obtain
that
\begin{eqnarray}
	-H_{n}'(t) & = &
	-2\exp\left(\frac{2}{3k}\,\Gkn(u_{n}(t))\right) 
	\langle\nabla\Gkn(u_{n}(t)),u_{n}'(t)\rangle
	\nonumber \\
	 & = & 2\exp\left(\frac{2}{3k}\,\Gkn(u_{n}(t))\right)
	 \left\|\nabla\Gkn(u_{n}(t))\right\|_{2}^{2/3}
	 \|u_{n}'(t)\|_{2}^{4/3}.
	\label{est:4/3}
\end{eqnarray}

Let us estimate the first two terms of this product for all
$t\in[0,\tsingn]$.  Let $J_{d,n}(t):=J_{d,n}(u_{n}(t))$ denote the
discrete jump heights of $u_{n}(t)$, defined according to
(\ref{defn:Jdn}).  From (\ref{est:f-vs-jh-min}) we have that
\begin{equation}
	\exp\left(\frac{2}{3k}\,\Gkn(u_{n}(t))\right)\geq
	\exp\left(\frac{2}{3}\log\left(\min_{d\in 
	D}|J_{d,n}(t)|\right)\right)=
	\left(\min_{d\in D}|J_{d,n}(t)|\right)^{2/3}.
	\label{est:hc1}
\end{equation}

Moreover, from (\ref{est:s-vs-jh-min2}) we have that
\begin{equation}
	\left\|\nabla\Gkn(u_{n}(t))\right\|_{2}^{2/3}\geq
	\left(\min_{d\in D}|J_{d,n}(t)|\right)^{-2/3}.
	\label{est:hc2}
\end{equation}

Plugging (\ref{est:hc1}) and (\ref{est:hc2}) into (\ref{est:4/3}) we 
obtain that
$$-H_{n}'(t)\geq 2\|u_{n}'(t)\|_{2}^{4/3}
\geq \|u_{n}'(t)\|_{2}^{4/3}
\quad\quad
\forall t\in[0,\tsingn].$$

Now we can conclude by integrating in $[s,t]$ and then applying
H\"{o}lder's inequality exactly as in the proof of
statement~(\ref{stat:u-hc}) of Proposition~\ref{prop:u}.\qed

\subsection{Well preparation}

In this section we prove our well preparation result
(Proposition~\ref{prop:wpl}).  To this end, in a time $S_{n}\to 0$ we
have three tasks to accomplish: extinguishing all vanishing jump
points, adjusting the energy, and remaining close enough to $u_{0}$.
In the next three lemmata we examine these three issues separately.
Then we make an alternate use of them in order to conclude the proof
of Proposition~\ref{prop:wpl}.

In the following, $J_{d}(u_{0})$ denotes the jump heights of $u_{0}$,
and $J_{d,n}(t):=J_{d,n}(u_{n}(t))$ denotes the discrete jump heights
of $u_{n}(t)$, defined according to (\ref{defn:Jdn}).

\begin{lemma}[Infinitesimal interval convergence]\label{lemma:uc}
	Let $D$, $D'$, $\{u_{0n}\}$, $u_{0}$, $u_{n}(t)$ be as in
	Theorem~\ref{thm:main}.  Let $\{S_{n}\}$ be any sequence of
	nonnegative real numbers such that $S_{n}\to 0$ as $n\to +\infty$.
	
	Then we have that
	\begin{equation}
		\lim_{n\to+\infty}\;\max_{t\in[0,S_{n}]}
		\|u_{n}(t)-u_{0}\|_{2}=0.
		\label{th:uc-sup}
	\end{equation}
	
	Moreover, for every $n$ large enough we have that $u_{n}(S_{n})\in
	PS_{D''(n),n}$ for some finite set $D''(n)$ (which may depend on
	$n$) such that $D'\subseteq D''(n)\subseteq D$.
\end{lemma}

\paragraph{\textmd{\emph{Proof}}} 

Let $w_{0n}\in\PCn$ denote the approximation of $u_{0}$ defined by
$$w_{0n}(x):=u_{0}\left(\frac{\lfloor nx\rfloor}{n}\right) 
\quad\quad
\forall x\in[0,1],\ \forall n\geq 1.$$

Since $u_{0}$ is already a piecewise constant function, we have that
$w_{0n}$ coincides with $u_{0}$ but for the subintervals corresponding
to elements of $D'$.  In other words, one can think of $w_{0n}$ as
obtained by moving every jump of $u$ in the point of the grid on its
left.  In particular, when $n$ is large enough we have that
\begin{equation}
	J_{d,n}(w_{0n})=J_{d}(u_{0})
	\quad\quad
	\forall d\in D'.
	\label{jwn=jdu0}
\end{equation}

Of course we have also that $w_{0n}\to u_{0}$ in $\LL$ as $n\to +\infty$.

Let $k':=|D'|$, let $K_{0}$ be the constant defined in
Lemma~\ref{lemma:jh-est}, and let us set 
$$c_{0}:=\min_{d\in D'}|J_{d}(u_{0})|,\quad c_{1}:=\max_{d\in
D'}|J_{d}(u_{0})|,\quad c_{2}:=\left(\min\left\{\frac{c_{0}}{12},
\frac{K_{0}}{12}\right\}\right)^{3},\quad
t_{0}:=\frac{c_{0}c_{2}}{8k'c_{1}}.$$

Let us consider the function
$$y_{n}(t):=\|u_{n}(t)-w_{0n}\|_{\LL((0,1))}^{2},$$
where the norm is intended with respect to the space variable, and let
\begin{equation}
	R_{n}:=\sup\left\{s\geq 0:y_{n}(t)\leq c_{2}
	\quad
	\forall t\in[0,s]\right\}.	
	\label{uc:defn-rn}
\end{equation}

Since
$$y_{n}(0)=\|u_{0n}-w_{0n}\|_{2}^{2}\leq
\left(\|u_{0n}-u_{0}\|_{2}+\|u_{0}-w_{0n}\|_{2}\strut\right)^{2},$$
it is easy to see that $y_{n}(0)\to 0$ as $n\to +\infty$, hence
$R_{n}$ is the supremum of a nonempty set when $n$ is large enough.
Let us fix $n_{0}$ big enough so that we can apply
Lemma~\ref{lemma:jh-est} for every $n\geq n_{0}$, and such that
$$y_{n}(0)\leq\frac{c_{2}}{2}
\quad\quad
\forall n\geq n_{0}.$$

We claim that for every $n\geq n_{0}$ we have that $R_{n}\geq t_{0}$, and
\begin{equation}
	y_{n}(t)\leq y_{n}(0)+\frac{4k'c_{1}}{c_{0}}t
	\quad\quad
	\forall t\in[0,t_{0}].
	\label{ead:est-yn}
\end{equation}

If we prove these claims, then (\ref{th:uc-sup}) is proved.  Indeed
$S_{n}\to 0$, hence $S_{n}\leq t_{0}$ for $n$ large enough, and
therefore 
$$\|u_{n}(t)-u_{0}\|_{2}\leq
\|u_{n}(t)-w_{0n}\|_{2}+\|w_{0n}-u_{0}\|_{2}\leq
\left(y_{n}(0)+\frac{4k'c_{1}}{c_{0}}S_{n}\right)^{1/2}+
\|w_{0n}-u_{0}\|_{2}$$
for every $t\in[0,S_{n}]$. Since all the terms in the 
right-hand side tend to zero, estimate (\ref{th:uc-sup}) is proved.

So we are left to prove these claims.  Let us consider
$t\in[0,R_{n}]$.  Let us apply (\ref{th:est-d-d}) with $v=u_{n}(t)$
and $w=w_{0n}$ (this can be done because the discontinuity sets of
$u(t)$ and $w_{0n}$ are contained in the fixed set $D$).  We obtain that
$$\min\left\{\strut K_{0},|J_{d,n}(t)-J_{d,n}(w_{0n})|\right\}\leq 
6\left[y_{n}(t)\right]^{1/3}\leq 6c_{2}^{1/3}
\quad\quad
\forall d\in D'.$$

The right-hand side is less than or equal to $K_{0}/2$ due to our
definition of $c_{2}$.  Combining with (\ref{jwn=jdu0}) it follows
that 
$$|J_{d,n}(t)-J_{d}(u_{0})|=
|J_{d,n}(t)-J_{d,n}(w_{0n})|\leq 6c_{2}^{1/3} 
\quad\quad
\forall d\in D'.$$

Therefore, from our definition of $c_{0}$ and $c_{2}$ it
follows that
\begin{equation}
	|J_{d,n}(t)|\geq|J_{d}(u_{0})|-|J_{d,n}(t)-J_{d}(u_{0})|\geq
	c_{0}-6c_{2}^{1/3}\geq c_{0}-\frac{c_{0}}{2}=\frac{c_{0}}{2}
	\label{ead:jump-big}
\end{equation}
for every $d\in D'$.  Now let us compute the time derivative of $y_{n}(t)$.
From equation (\ref{pbm:main-eq}) and formula (\ref{defn:grad-gkn}) we
obtain that
\begin{eqnarray*}
	y_n'(t) &=& 2\int_0^1 \left(u_n( x, t)-w_{0n}(x)
	\right) \frac {d}{dt} u_n(x,t)\, dx\\
	& = & 2\int_0^1 \left(u_n( x, t)-w_{0n}(x) \right)\cdot nD^{-1/n}\left[
	\frac{D^{1/n}u_{n}(x,t)}{1+|D^{1/n}u_{n}(x,t)|^{2}}\right]\,dx.
\end{eqnarray*}

Now we apply the discrete version of the integration-by-parts formula 
(which actually is a simple algebraic manipulation of finite sums). 
We do not have boundary terms because $D^{1/n}u_{n}(x,t)$ is zero 
both in the last subinterval $(1-1/n,1)$ and in $(-1/n,0)$.
We obtain that
\begin{eqnarray*}
	y_n'(t) & = & -2\int_0^1 \left(
	D^{1/n}u_n(x,t)-D^{1/n}w_{0n}(x) \right)\cdot n\,
	\frac{D^{1/n}u_{n}(x,t)}{1+|D^{1/n}u_{n}(x,t)|^{2}}\,dx  \\
	\noalign{\vspace{1ex}}
	 & \leq & 2\int_{0}^{1}|D^{1/n} w_{0n}(x)|\cdot n\,
	\frac{|D^{1/n}u_{n}(x,t)|}{1+|D^{1/n}u_{n}(x,t)|^{2}}\,dx  \\
	\noalign{\vspace{1ex}}
	 & = & 2\sum_{d\in D'} |J_{d,n}(w_{0n})|\cdot
	 \frac{n^2|J_{d,n}(t)|}{1+n^2|J_{d,n}(t)|^2}  \\
	 &\leq& 2\sum_{d\in D'} |J_{d,n}(w_{0n})|\cdot
	 \frac{1}{|J_{d,n}(t)|}.
\end{eqnarray*}

The terms of the last sum can be estimated using (\ref{jwn=jdu0}), our
definition of $c_{1}$, and (\ref{ead:jump-big}).  We obtain that
$y_{n}'(t)\leq 4k'c_{1}/c_{0}$ for every $t\in[0,R_{n}]$, hence the
estimate in (\ref{ead:est-yn}) holds true for every $t\in[0,R_{n}]$.

Let us assume now that $R_{n}<t_{0}$ for some $n\geq n_{0}$. Due to 
the maximality of $R_{n}$ we obtain that
$$c_{2}=y_{n}(R_{n})\leq y_{n}(0)+\frac{4k'c_{1}}{c_{0}}R_{n}<
\frac{c_{2}}{2}+\frac{4k'c_{1}}{c_{0}}t_{0}=c_{2},$$
which is a contradiction. This completes the proof of our claims, 
hence also of (\ref{th:uc-sup}).

It remains to understand the location of discrete jump points of
$u_{n}(S_{n})$.  First of all, the creation of new jump points is
forbidden by statement~(\ref{stat:monotonicity}) of
Theorem~\ref{thmbibl:semidiscrete}, and this proves that
$D''(n)\subseteq D$.  On the other hand, we have that $u_{n}(S_{n})\to
u_{0}$, hence statement~(\ref{stat:jc}) of Lemma~\ref{lemma:unif-est}
implies that $D'\subseteq D''(n)$ when $n$ is large enough.\qed

\begin{lemma}[Energy adjustment]\label{lemma:ead}
	Let $D$, $D'$, $\{u_{0n}\}$, $u_{0}$, $u_{n}(t)$, be as in
	Theorem~\ref{thm:main}, and let $k:=|D|$. Let $\{S_{n}\}$ be any 
	sequence such that $S_{n}\to 0$ and
	\begin{equation}
		\lim_{n\to +\infty}ne^{-2nS_{n}}=0.
		\label{hp:sn-lim}
	\end{equation}
	Let $\tsingn$ be defined as in (\ref{defn:tsingn}), and let us assume that
	\begin{equation}
		\tsingn>S_{n}
		\quad
		\mbox{for every $n$ large enough}.
		\label{hp:sn-tsingn}
	\end{equation}
	
	Then we have that
	\begin{equation}
		\lim_{n\to+\infty}\Gkn(u_{n}(S_n)) =\Gki(u_{0}),
		\label{th:ead-gkn}
	\end{equation}
	where of course $\Gki(u_{0})=-\infty$ if $D'$ is strictly 
	contained in $D$.
\end{lemma}

\paragraph{\textmd{\emph{Proof}}} 

Since $S_{n}\to 0$, from Lemma~\ref{lemma:uc} we know that
$u_{n}(S_{n})\to u_{0}$.  Therefore, the Gamma-convergence of $\Gkn$
to $\Gki$ implies that
$$\liminf_{n\to +\infty}\Gkn(u_{n}(S_{n}))\geq\Gki(u_{0}).$$

So we are left to prove the opposite inequality with the $\limsup$, 
which in turn is equivalent to show that
\begin{equation}
	\limsup_{n\to +\infty}\Gkn(u_{n}(S_{n}))\leq M
	\label{ead:limsup}
\end{equation}
for every $M>\Gki(u_{0})$.

Let us fix any such $M$. To begin with, we claim that
\begin{equation}
	\frac{1}{2} \sum_{d\in D}
	\log\left(\frac{1}{n^{2}}+J_{d,n}^{2}(t)\right) \leq M
	\quad\quad
	\forall t\in[0,S_{n}]
	\label{ead:old-rn}
\end{equation}
for every $n$ large enough.  Indeed let us assume that this is not the
case.  Then there exists a sequence $\{t_{n}\}$, with
$t_{n}\in[0,S_{n}]$ for every $n\geq 1$, such that
\begin{equation}
	\frac{1}{2} \sum_{d\in D}
	\log\left(\frac{1}{n^{2}}+J_{d,n}^{2}(t_{n})\right) > M
	\quad\quad
	\mbox{for infinitely many $n$'s}.
	\label{ead:absurd}
\end{equation}

On the other hand, from Lemma~\ref{lemma:uc} we deduce that 
$u_{n}(t_{n})\to u_{0}$, hence from the jump convergence 
(\ref{th:jump-conv}) we obtain that $J_{d,n}(t_{n})\to J_{d}(u_{0})$ 
for every $d\in D$, where of course $J_{d}(u_{0})=0$ if $d\in 
D\setminus D'$. In particular we have that
$$\lim_{n\to +\infty}\frac{1}{2} \sum_{d\in D}
\log\left(\frac{1}{n^{2}}+J_{d,n}^{2}(t_{n})\right)= \Gki(u_{0})< M,$$
which contradicts (\ref{ead:absurd}).

From now on we work in the interval $[0,S_{n}]$, and we assume $n$ to
be large enough so that (\ref{ead:old-rn}) and (\ref{hp:sn-tsingn})
hold true.  Since $S_{n}<\tsingn$, in this interval we know that
$u_{n}(t)\in\PSDn$, hence we can use the estimates of
Lemma~\ref{lemma:fund-est}.  Let $SQ_{n}(t):=SQ_{n}(u_{n}(t))$ denote
the subcritical incremental quotient of $u_{n}(t)$, defined according
to (\ref{defn:SQn}).  Let us estimate the time derivative of the
function $t\to\Gkn(u_{n}(t))$.  From equation (\ref{pbm:main-eq-gkn})
and estimate (\ref{est:slope-vs-sq}) we have that
\begin{equation}
	\frac{d}{dt}\Gkn(u_{n}(t))\ =\ 
	-\left\|\nabla\Gkn(u_{n}(t))\right\|_{2}^{2}
	\ \leq\ -n^{2}[SQ_{n}(t)]^{2}
	\label{est:ead1}
\end{equation}
for every $t\in[0,S_{n}]$.  On the other hand, from
(\ref{est:f-above}) and (\ref{ead:old-rn}) we have that
\begin{eqnarray}
	\Gkn(u_{n}(t)) & \leq & \frac{n}{2}\log\left(1+[SQ_{n}(t)]^{2}\right)+
	\frac{1}{2} \sum_{d\in
	D}\log\left(\frac{1}{n^{2}}+J_{d,n}^{2}(t)\right)  
	\nonumber\\
	 & \leq & 
	 \frac{n}{2}\log\left(1+[SQ_{n}(t)]^{2}\right)+M,
	 \label{est:gkn}
\end{eqnarray}
hence
$$[SQ_{n}(t)]^{2}\ \geq\ 
\exp\left(\frac{2}{n}\left[\Gkn(u_{n}(t))-M\right]\right)-1\ \geq\ 
\frac{2}{n}\left[\Gkn(u_{n}(t))-M\right].$$

Plugging this estimate into (\ref{est:ead1}) we obtain that
$$\frac{d}{dt}\Gkn(u_{n}(t))\leq
-2n\left[\Gkn(u_{n}(t))-M\right]
\quad\quad
\forall t\in[0,S_{n}],$$
hence
$$\Gkn(u_{n}(t))\leq\left[\Gkn(u_{0n})-M\right]e^{-2nt}+M
\quad\quad
\forall t\in[0,S_{n}].$$

Let us estimate $\Gkn(u_{0n})$.  From (\ref{est:gkn}) with $t=0$ we
have that
$$\Gkn(u_{0n})\
\leq\ \frac{n}{2}\log\left(1+[SQ_{n}(0)]^{2}\right)+M 
\ \leq\ \frac{n}{2}\log 2+M,$$
hence
$$\Gkn(u_{n}(S_{n}))\leq
\frac{\log 2}{2} \cdot n e^{-2nS_{n}}+M.$$

If $S_{n}$ satisfies (\ref{hp:sn-lim}), then the first term in the
right-hand side tends to 0.  This completes the proof of
(\ref{ead:limsup}), hence also the proof of (\ref{th:ead-gkn}).\qed

\begin{lemma}[Discrete jump extinction]\label{lemma:dje}
	Let $D$, $D'$, $\{u_{0n}\}$, $u_{0}$, $u_{n}(t)$ be as in
	Theorem~\ref{thm:main}.  Let us assume that $D'$ is strictly
	contained in $D$, and let $\tsingn$ be the first time when a
	discrete jump disappears, defined according to
	(\ref{defn:tsingn}).

	Then we have that
	\begin{equation}
		\lim_{n\to +\infty}\tsingn=0.
		\label{th:dje-tsingn}
	\end{equation}
\end{lemma}

\paragraph{\textmd{\emph{Proof}}} 

Let us assume by contradiction that (\ref{th:dje-tsingn}) is false. 
This is equivalent to say that there exist $\delta>0$ and a 
subsequence (not relabeled) such that
\begin{equation}
	\tsingn\geq\delta>0
	\quad
	\mbox{for every $n$ large enough}.
	\label{dje:absurd}
\end{equation}

Let us take any sequence $S_{n}\to 0$ satisfying (\ref{hp:sn-lim}), for 
example $S_{n}=n^{-1/2}$. Due to (\ref{dje:absurd}) this sequence 
satisfies also (\ref{hp:sn-tsingn}), hence we can apply 
Lemma~\ref{lemma:ead} to this sequence.  Since $D'$ is
strictly contained in $D$, we have that $\Gki(u_{0})=-\infty$, hence
(\ref{th:ead-gkn}) reads as
$$\lim_{n\to+\infty}\Gkn(u_{n}(S_{n}))=-\infty.$$

Moreover, from (\ref{hp:sn-tsingn}) we have that
$u_{n}(S_{n})\in\PSDn$ for every $n$ large enough.  Therefore, up to
replacing the initial sequence $u_{0n}\to u_{0}$ with the sequence
$u_{n}(S_{n}) \to u_{0}$ (the convergence to $u_{0}$ is due to 
Lemma~\ref{lemma:uc}), we can always assume that the sequence of
initial data satisfies (\ref{hp:u0n-in-psdn}) and
\begin{equation}
	\lim_{n\to+\infty}\Gkn(u_{0n})=-\infty.
	\label{hp:dje-gkn}
\end{equation}

Thus from now on we work under this assumption.  Let us estimate the
time derivative of the function $t\to \Gkn(u_{n}(t))$ in the interval 
$[0,\tsingn]$.  Using equation
(\ref{pbm:main-eq-gkn}), and estimate (\ref{est:s-vs-jh-min2}), we
find that
$$\frac{d}{dt}\Gkn(u_{n}(t))= -\left\|\nabla
\Gkn(u_{n}(t))\right\|_{2}^{2}\leq 
-\left(\min_{d\in D}|J_{d,n}(t)|\right)^{-2}.$$

Combining with (\ref{est:f-vs-jh-min}) we deduce that
$$\frac{d}{dt}\Gkn(u_{n}(t))\leq
-\exp\left(-\frac{2}{k}\,\Gkn(u_{n}(t))\right) .$$

Integrating this differential inequality in $[0,\tsingn]$ we
obtain that
$$\frac 2 k\,\tsingn \leq \exp\left(\frac 2 k\, \Gkn(u_{n}(0))\right)- 
\exp\left(\frac 2 k\, \Gkn(u_{n}(\tsingn))\right) \leq 
\exp\left(\frac 2 k\, \Gkn(u_{0n})\right).$$

Thanks to assumption (\ref{hp:dje-gkn}), this contradicts 
(\ref{dje:absurd}), hence it proves 
(\ref{th:dje-tsingn}).\qed

\subsubsection*{Proof of Proposition~\ref{prop:wpl}}

We argue by induction on $k-k'$, where $k:=|D|$ and $k':=|D'|$.

Let us assume that $k-k'=0$, namely $D=D'$.  In this case we claim
that conclusions (\ref{th:wpl-dje}) through (\ref{th:wpl-sup}) hold
true for every sequence $S_{n}\to 0$ satisfying (\ref{hp:sn-lim}), for
example $S_{n}:=n^{-1/2}$.  First of all, from Lemma~\ref{lemma:uc} we
obtain that (\ref{th:wpl-sup}) holds true, and $u_{n}(S_{n})\in
PS_{D',n}$ when $n$ is large enough (because $D=D'$).  This proves
(\ref{th:wpl-dje}), and implies that assumption (\ref{hp:sn-tsingn})
is satisfied.  Therefore, we can apply Lemma~\ref{lemma:ead} and
deduce (\ref{th:ead-gkn}).  Since $k=k'$, this proves
(\ref{th:wpl-wp}).

If $k-k'>0$, then we begin by applying Lemma~\ref{lemma:dje}. We 
obtain the extinction of at least one 
discrete jump in a time $\tsingn\to 0$. Thus from 
Lemma~\ref{lemma:uc} we have also that
\begin{equation}
   u_{n}(\tsingn)\to u_{0}
   \label{wpl:tsingn-to-u0}
\end{equation}
and every element of this new sequence belongs to $PS_{D''(n),n}$ for
some $D''(n)$ strictly contained in $D$ and such that $D''(n)\supseteq
D'$.  Up to subsequences we can assume that $D''(n)=:D''$ is
independent of $n$.  Since the number of possible choices of $D''$ is
finite, we have only finitely many subsequences to consider.
Therefore, it is enough to conclude on all such subsequences.

Setting $k'':=|D''|$, all these subsequences of (\ref{wpl:tsingn-to-u0})
satisfy the same assumptions of the initial sequence with $D''$
instead of $D$, and in particular with $k''-k'<k-k'$.  Therefore, the
conclusion follows from the inductive assumption.\qed

\subsection{Convergence up to collisions}

\subsubsection*{Proof of Proposition~\ref{prop:fino-t1}}

\paragraph{\textmd{\emph{``Safe intervals''}}}

Let $J_{d}(u_{0})$ be the jump heights of $u_{0}$,
and let $J_{d,n}(t):=J_{d,n}(u_{n}(t))$ be the
discrete jump heights of $u_{n}(t)$, defined according to 
(\ref{defn:Jdn}).  

We say that $[0,T_{0}]$ (with $T_{0}>0$) is a safe interval if there
exists a positive real number $c_{0}$, and a positive integer $n_{0}$
(both may depend on $T_{0}$) such that
\begin{equation}
	|J_{d,n}(t)|\geq c_{0}
	\quad\quad
	\forall d\in D,\ \forall t\in[0,T_{0}],\ \forall n\geq n_{0}.
	\label{defn:safe-int}
\end{equation}

We claim that safe intervals do exist.  Indeed let us assume that this
is not the case.  Then there exist a sequence $\{d_{k}\}\subseteq D$,
a sequence $n_{k}\to +\infty$ of positive integers, and a sequence
$t_{k}\to 0$ of positive times such that
\begin{equation}
	\lim_{k\to +\infty}J_{d_{k},n_{k}}(t_{k})=0.
	\label{absurd:si}
\end{equation}

Up to subsequences, we can assume that $d_{k}$ does not depend on $k$.
Since $t_{k}\to 0$, Lemma~\ref{lemma:uc} implies that
$u_{n_{k}}(t_{k})\to u_{0}$.  At this point the jump convergence
(\ref{th:jump-conv}) contradicts (\ref{absurd:si}). From 
(\ref{defn:safe-int}) it is also clear that
\begin{equation}
	\tsingn\geq T_{0}
	\quad\quad
	\mbox{for every $n$ large enough}.
	\label{est:t1-tsingn-t0}
\end{equation}

\paragraph{\textmd{\emph{Boundedness and compactness in safe intervals}}}

Let $[0,T_{0}]$ be a safe interval according to 
(\ref{defn:safe-int}). In this part of the proof we show that there 
exist real constants $c_{1}$ and $c_{2}$ such that
\begin{equation}
	c_{1}\leq\Gkn(u_{n}(t))\leq c_{2}
	\quad\quad
	\forall t\in[0,T_{0}],\ \forall n\geq 1,
	\label{est:gkn-bounded}
\end{equation}
and that there exists $v\in C^{0}\left([0,T_{0}];\LL\right)$ such that
(up to subsequences, which we do not relabel)
\begin{equation}
	u_{n}\to v
	\quad\mbox{in }C^{0}\left([0,T_{0}];\LL\right).
	\label{fino-t1:un-to-v}
\end{equation}

Indeed from the monotonicity of the function $t\to\Gkn(u_{n}(t))$ we
have that
$$\Gkn(u_{n}(t))\leq\Gkn(u_{0n})
\quad\quad
\forall t\geq 0,$$
and the right-hand side is bounded from above because of
(\ref{hp:main-gkn}).  Moreover, (\ref{est:f-vs-jh-min}) and the safe
interval assumption (\ref{defn:safe-int}) imply that $\Gkn(u_{n}(t))
\geq k\log c_{0}$ 
for every $n\geq n_{0}$.  This completes the proof of
(\ref{est:gkn-bounded}).  Now we exploit a compactness argument.
\begin{itemize}
	\item For every $t\in[0,T_{0}]$ (and actually for every $t\geq 0$)
	the sequence $\{u_{n}(t)\}$ is relatively compact in $\LL$.  Indeed
	from statement~(\ref{stat:unif-est}) of 
	Lemma~\ref{lemma:unif-est} (applied with $D'=D$) we have that
	$$\sup_{n\geq
	1}\left(\|u_{0n}\|_{\infty}+
	\left\|D^{1/n}u_{0n}\right\|_{1}\right)
	<+\infty,$$	
	hence from statements~(3) and~(4) of
	Theorem~\ref{thmbibl:semidiscrete} we deduce that 
	$$\sup_{n\geq
	1}\left(\|u_{n}(t)\|_{\infty}+
	\left\|D^{1/n}u_{n}(t)\right\|_{1}\right)
	<+\infty.$$
	
	In other words, we control the $L^{\infty}$-norm and the total
	variation of $u_{n}(t)$ (as functions of the space variable).
	This guarantees the required compactness.

	\item The functions $u_{n}:[0,T_{0}]\to\LL$ are 1/2-H\"{o}lder
	continuous, with equi-bounded H\"{o}lder constants, because of
	(\ref{est:gkn-bounded}) and statement~(\ref{stat:holder-1/2}) of
	Theorem~\ref{thmbibl:semidiscrete}.  Alternatively, they are
	1/4-H\"{o}lder continuous, with equi-bounded H\"{o}lder constants,
	because of (\ref{est:t1-tsingn-t0}), (\ref{hp:main-gkn}), and 
	Proposition~\ref{prop:hc}.
\end{itemize}

Therefore, Ascoli's Theorem implies that the sequence $\{u_{n}(t)\}$ 
is relatively compact in $C^{0}\left([0,T_{0}];\LL\right)$. This 
proves (\ref{fino-t1:un-to-v}).

\paragraph{\textmd{\emph{Passing to the limit in safe intervals}}}

Let $[0,T_{0}]$ be a safe interval according to (\ref{defn:safe-int}),
and let $v(t)$ be any limit point of the sequence $u_{n}(t)$.  In this
part of the proof we show that $v(t)=u(t)$ in the safe interval
$[0,T_{0}]$.  As a consequence, we obtain also that
(\ref{fino-t1:un-to-v}) holds true for the whole sequence, and not
only up to subsequences.

To this end, we write the differential equations in integral form, as
in the theory of maximal slope curves.  From statement~(1) of
Proposition~\ref{prop:cmp}, we know that equation
(\ref{pbm:main-eq-gkn}) implies that (and actually is equivalent to)
\begin{equation}
	\Gkn(u_n(s))-\Gkn(u_n(t))= \frac 1 2 \int_s^t
	\|u_n'(\tau)\|_{2}^2\,d\tau+ \frac 1 2 \int_s^t \left\|\nabla \Gkn
	(u_n(\tau))\right\|_{2}^2\,d\tau
	\label{eqn:cmp}
\end{equation}
for every $0\leq s\leq t\leq T_{0}$.

From (\ref{est:gkn-bounded}) we know that the left-hand side of
(\ref{eqn:cmp}) is bounded from above.
In particular, setting $s=0$ and $t=T_{0}$, we obtain that
\begin{equation}
	\sup_{n\geq 1}\int_{0}^{T_{0}}\|u_n'(\tau)\|_{2}^2\,d\tau<+\infty,
	\label{est:t1-deriv}
\end{equation}
\begin{equation}
	\sup_{n\geq 1}\int_{0}^{T_{0}}\left\|\nabla \Gkn
	(u_n(\tau))\right\|_{2}^2\,d\tau<+\infty.
	\label{est:t1-slope}
\end{equation}

From (\ref{est:t1-deriv}) we easily deduce that $v\in 
H^{1}\left((0,T_{0});\LL\right)$, and
\begin{equation}
	\liminf_{n\to +\infty}\int_s^t \|u_n'(\tau)\|_{2}^2\,d\tau\geq
	\int_s^t \|v'(\tau)\|_{2}^2\,d\tau
	\label{est:t1-liminf-deriv}
\end{equation}
for every $0\leq s\leq t\leq T_{0}$.

Let us consider now the second term in the right-hand side of
(\ref{eqn:cmp}).  Due to the safe interval assumption
(\ref{defn:safe-int}) and (\ref{est:t1-tsingn-t0}), we can apply Proposition~\ref{prop:cmp-limit}.  From
(\ref{th:slope-liminf}) we deduce that 
$$\liminf_{n\to +\infty}\left\|\nabla\Gkn(u_{n}(t))\right\|_{2}\geq
\left\|\nabla\Gki(v(t))\right\|_{2} \quad\quad \forall
t\in[0,T_{0}].$$

Thus from Fatou's Lemma it follows that
\begin{eqnarray}
	\liminf_{n\to +\infty}\int_s^t \left\|\nabla \Gkn
		(u_n(\tau))\right\|_{2}^2\,d\tau & \geq & \int_s^t
		\left(\liminf_{n\to +\infty}\left\|\nabla \Gkn
		(u_n(\tau))\right\|_{2}^2\right)\,d\tau
	\nonumber  \\
	 & \geq & \int_s^t \left\|\nabla \Gki
	(v(\tau))\right\|_{2}^2\,d\tau
	\label{est:t1-liminf-slope}
\end{eqnarray}
for every $0\leq s\leq t\leq T_{0}$.

Now we consider the left-hand side of (\ref{eqn:cmp}).  The functions
$t\to\Gkn(u_{n}(t))$ are equi-bounded and nonincreasing.  By the usual
compactness result for monotone functions (known as Helly's Lemma, see
\cite[Lemma~3.3.3]{AGS} ) there exists a nonincreasing function
$\psi:[0,T_{0}]\to\re$ such that (up to subsequences)
\begin{equation}
	\lim_{n\to +\infty}\Gkn(u_{n}(t))=\psi(t)
	\quad\quad
	\forall t\in[0,T_{0}].
	\label{defn:psi}
\end{equation}

Now we can take the liminf of both sides of (\ref{eqn:cmp}).  Thanks
to (\ref{est:t1-liminf-deriv}), (\ref{est:t1-liminf-slope}), and
(\ref{defn:psi}) we obtain that
\begin{equation}
	\psi(s)-\psi(t)\geq \frac 1 2
	\int_s^t \|v'(\tau)\|_{2}^2\,d\tau+ 
	\frac 1 2 \int_s^t \left\|\nabla \Gki
	(v(\tau))\right\|_{2}^2\,d\tau.
	\label{eqn:cmp-lim-psi}
\end{equation}

It remains to characterize the function $\psi(t)$.  Coming back to
(\ref{est:t1-slope}), and exploiting once more Fatou's Lemma, we
obtain that 
$$\int_{0}^{T_{0}} \left(\liminf_{n\to
+\infty}\left\|\nabla \Gkn (u_n(\tau))\right\|_{2}^2\right)\,d\tau\leq
\liminf_{n\to +\infty}\int_{0}^{T_{0}} \left\|\nabla \Gkn
(u_n(\tau))\right\|_{2}^2\,d\tau<+\infty.$$

Therefore there exists a set $E\subseteq[0,T_{0}]$, with Lebesgue 
measure equal to 0, such that
$$\liminf_{n\to +\infty}\left\|\nabla \Gkn (u_n(t))\right\|_{2}<+\infty
\quad\quad
\forall t\in[0,T_{0}]\setminus E.$$

As a consequence, for every $t\in[0,T_{0}]\setminus E$ there exists a 
($t$-dependent) sequence $n_{h}\to +\infty$ such that
$$\sup_{h\in\n}\left\|\nabla
G^{(k)}_{n_{h}}(u_{n_{h}}(t))\right\|_{2}<+\infty.$$

On this subsequence we can apply statement~(2) of 
Proposition~\ref{prop:cmp-limit} and deduce that
$$v(t)=\lim_{n\to +\infty}u_{n}(t)=
\lim_{h\to +\infty}u_{n_{h}}(t)\in\PCD,$$
$J_{d}(v(t))\geq c_{0}$ for every $d\in D$, and
$$\psi(t)=\lim_{n\to +\infty}\Gkn(u_{n}(t))=
\lim_{h\to +\infty}G^{(k)}_{n_{h}}(u_{n_{h}}(t))=
\Gki(v(t)).$$

We have thus proved that
\begin{equation}
	\psi(t)=\Gki(v(t))
	\quad\quad
	\forall t\in[0,T_{0}]\setminus E.
	\label{eqn:psi-ae}
\end{equation}

The last step is to prove the same equality for every $t\in(0,T_{0})$.
To this end we remark that $v(t)$ is a continuous function with values
in $\PCD$, and $\Gki$ is continuous in $\PCD$.  It follows that the
right-hand side of (\ref{eqn:psi-ae}) is a continuous function.
Therefore, in (\ref{eqn:psi-ae}) we have a continuous function and a
monotone function which coincide almost everywhere in $[0,T_{0}]$,
hence they coincide everywhere in $(0,T_{0})$.

Coming back to (\ref{eqn:cmp-lim-psi}), we have proved that
$$\Gki(v(s))-\Gki(v(t))\geq \frac 1 2 \int_s^t \|v'(\tau)\|_{2}^2\,d\tau+
\frac 1 2 \int_s^t \left\|\nabla \Gki (v(\tau))\right\|_{2}^2\,d\tau$$
for every $0< s\leq t< T_{0}$.  From statement~(2) of
Proposition~\ref{prop:cmp}, this is equivalent to say that $v(t)$
coincides with $u(t)$ in $[0,T_{0}]$.

We have also proved the energy convergence (\ref{th:t1-gkn<}) for 
every $t\in(0,T_{0})$.

\paragraph{\textmd{\emph{Continuation up to first jump extinction}}}

Let $T_{0\infty}$ be the supremum of all $T_{0}>0$ such that 
$[0,T_{0}]$ is a safe interval according to (\ref{defn:safe-int}). 
From (\ref{est:t1-tsingn-t0}), and the convergence results on safe 
intervals, it is easy to see that
\begin{eqnarray}
	 & \displaystyle{\liminf_{n\to +\infty}\tsingn\geq T_{0\infty},} & 
	\label{est:t1-liminf-t0i}  \\
	\noalign{\vspace{1ex}}
	 & \displaystyle{\lim_{n\to +\infty}u_{n}(t)=u(t)\in\PCD}
	\quad\quad
	\forall t\in[0,T_{0\infty}), & 
	\label{est:t1-lim-u-t0i}  \\
	\noalign{\vspace{1ex}}
	 & \displaystyle{\lim_{n\to +\infty}\Gkn(u_{n}(t))=\Gki(u(t))}
	\quad\quad
	\forall t\in(0,T_{0\infty}). & 
	\label{est:t1-lim-G-t0i}
\end{eqnarray}

Let $\{R_{m}\}$ be an increasing sequence of positive real numbers 
such that $R_{m}\to T_{0\infty}$ as $m\to +\infty$, and let us set
$$A_{m,n}:=\max\left\{\strut\|u_{n}(t)-u(t)\|_{2}:0\leq 
t\leq\min\{R_{m},(1-1/n)\tsingn\}\right\}.$$

Since $[0,R_{m}]$ is a safe interval for each $m$, we have that
$$A_{m,n}\stackrel{n\to+\infty}{\longrightarrow}0
\stackrel{m\to+\infty}{\longrightarrow}0.$$

Therefore, Lemma~\ref{lemma:amn} (standard conclusion) implies the
existence of a sequence $m_{n}\to +\infty$ of positive integers such
that $A_{m_{n},n}\to 0$ as $n\to +\infty$.  We claim that
$$T_{n}:=\min\left\{R_{m_{n}},\left(1-\frac{1}{n}\right)\tsingn\right\}$$
is a sequence which satisfies (\ref{th:fino-t1-tn}) through
(\ref{th:fino-t1-sup}).  Indeed (\ref{th:fino-t1-tn}) is trivial, and
(\ref{th:fino-t1-sup}) is equivalent to say that $A_{m_{n},n}\to 0$.
It remains to prove (\ref{th:fino-t1-lim}).  From
(\ref{est:t1-liminf-t0i}) we easily deduce that $T_{n}\to
T_{0\infty}$.  Therefore, proving (\ref{th:fino-t1-lim}) is equivalent
to show that $T_{0\infty}=\tsing$.

To this end, from (\ref{est:t1-lim-u-t0i}) we immediately deduce that
$T_{0\infty}\leq\tsing$. Moreover, since
$$ \| u_{n}(T_{n})- u(T_{0\infty})\|_{2} \leq \| u_{n}(T_{n})- 
u(T_{n})\|_{2} +\| u(T_{n})- u(T_{0\infty})\|_{2},$$
from (\ref{th:fino-t1-sup}) and the continuity of $u$ we 
have in particular that
\begin{equation}
	u_{n}(T_{n})\to u(T_{0\infty}).
	\label{eqn:t1-tn-un}
\end{equation}

Let us assume now by contradiction that $T_{0\infty}<\tsing$, hence 
$u(T_{0\infty})\in\PCD$. In this case we claim that there exists 
$\delta>0$ such that $[0,T_{0\infty}+\delta]$ is a safe interval, and 
this contradicts the maximality of $T_{0\infty}$.

In order to prove the claim, we argue as in the first paragraph of 
the proof. If the claim is false, then there exist a sequence 
$\{d_{k}\}\subseteq D$, a sequence $n_{k}\to +\infty$ of positive 
integers, and a sequence $\{t_{k}\}$ of times such that $t_{k}\in[0, 
T_{0\infty}+1/k]$ for every $k\geq 1$, and such that 
(\ref{absurd:si}) holds true. Up to subsequences, we can assume that 
$d_{k}$ does not depend on $k$, and that $t_{k}$ tends to some limit 
$t_{\infty}\in[0,T_{0\infty}]$. We can also assume that either 
$t_{k}\in[0,T_{n_{k}}]$ for every $k\geq 1$, or  
$t_{k}\in[T_{n_{k}},T_{0\infty}+1/k]$ for every $k\geq 1$.

In the first case we have that $u_{n_{k}}(t_{k})\to u(t_{\infty})$
because of (\ref{th:fino-t1-sup}). In the second case we have that
$u_{n_{k}}(t_{k})\to u(t_{\infty})=u(T_{0\infty})$ because of
Lemma~\ref{lemma:uc} applied to the sequence of ``initial data''
(\ref{eqn:t1-tn-un}).  In both cases the limit lies in $\PCD$, hence the
jump convergence (\ref{th:jump-conv}) contradicts (\ref{absurd:si}).

Therefore, we have proved that
$T_{0\infty}=\tsingn$, hence also (\ref{th:fino-t1-lim}).  At this
point conclusion (\ref{th:t1-gkn<}) is exactly
(\ref{est:t1-lim-G-t0i}).\qed

\subsection{Proof of main results}

\subsubsection*{Proof of Theorem~\ref{thm:main}}

\subparagraph{\textmd{\emph{Global-in-time $\LL$-convergence}}}

We argue by induction on $k'=|D'|$.  If $k'=0$, then $u_{0}$ is a constant
function, and $u(t)\equiv u_{0}$ is the stationary solution.  In order
to prove (\ref{th:main}) it is therefore enough to show that the
function $t\to \|u_n(x,t)-u_{0}\|_{2}$ is nonincreasing.  This is true
because in this case $u_{n}(t)-u_{0}$ is once again a solution of
equation (\ref{pbm:main-eq}), hence its $L^{2}$-norm is a
nonincreasing function of time because of statement~(3) of
Theorem~\ref{thmbibl:semidiscrete}.

Now let us consider the case where $k'>0$. In this case we argue in 
three steps.

First of all we exploit Proposition~\ref{prop:wpl} in order to ``well
prepare'' the sequence of initial data.  Let $\{S_{n}\}$ be the
sequence of times provided by Proposition~\ref{prop:wpl}.  We have that
$u_{n}(S_{n})\to u_{0}$ is a ``well prepared'' sequence with respect
to the discontinuity set $D'$, in the sense that the elements of
this sequence lie in the corresponding space $PS_{D',n}$ when $n$ is
large enough, and their $k'$-energies converge to the $k'$-energy of
$u_{0}$. Now let us observe that
$$\max_{t\in[0,S_{n}]}\|u_{n}(t)-u(t)\|_{2}\leq
\max_{t\in[0,S_{n}]}\|u_{n}(t)-u_{0}\|_{2}+
\max_{t\in[0,S_{n}]}\|u_{0}-u(t)\|_{2}.$$

The first term in the right-hand side tends to 0 as $n\to +\infty$ 
because of (\ref{th:wpl-sup}). The second term tends to 0 because 
$S_{n}\to 0$ and $u$ is continuous. It follows that
\begin{equation}
	\lim_{n\to +\infty}\,
	\max_{t\in[0,S_{n}]}\|u_{n}(t)-u(t)\|_{2}=0.
	\label{main:0-sn}
\end{equation}

The second step is to apply Proposition~\ref{prop:fino-t1} to the
``well prepared'' sequence of ``initial data'' $u_{n}(S_{n})\to
u_{0}$.  Let $\{T_{n}\}$ be the
sequence of times provided by Proposition~\ref{prop:fino-t1}. Then 
(\ref{th:fino-t1-sup}) reads as
\begin{equation}
	\lim_{n\to +\infty}\,
	\max_{t\in[0,T_{n}]}\|u_{n}(S_{n}+t)-u(t)\|_{2}=0.
	\label{main:fino-t1}
\end{equation}

On the other hand we have that
$$\max_{t\in[S_{n},S_{n}+T_{n}]}\|u_{n}(t)-u(t)\|_{2}\leq
\max_{t\in[0,T_{n}]}\|u_{n}(S_{n}+t)-u(t)\|_{2}+
\max_{t\in[0,T_{n}]}\|u(t)-u(S_{n}+t)\|_{2}.$$

The first term in the right-hand side tends to 0 as $n\to +\infty$ 
because of (\ref{main:fino-t1}). The second term tends to 0 because 
$S_{n}\to 0$ and $u$ is uniformly continuous. It follows that
\begin{equation}
	\lim_{n\to +\infty}\;
	\max_{t\in[S_{n},S_{n}+T_{n}]}\|u_{n}(t)-u(t)\|_{2}=0.
	\label{main:sn-tn}
\end{equation}

From (\ref{th:fino-t1-lim}) we have also that $S_{n}+T_{n}\to\tsing$. 
Therefore, since
$$\|u_{n}\left(S_{n}+T_{n}\right)- u(\tsing)\|_{2}\leq
\|u_{n}\left(S_{n}+T_{n}\right)- u(S_{n}+T_{n})\|_{2} + 
\|u\left(S_{n}+T_{n}\right)- u(\tsing)\|_{2},
$$
from (\ref{main:sn-tn}) and the continuity of $u$ it follows that
\begin{equation}
	u_{n}\left(S_{n}+T_{n}\right)\to u(\tsing).
	\label{main:us-sn+tn}
\end{equation}

In the last step we consider this sequence.  Due to
(\ref{th:fino-t1-tn}), all the elements of this sequence belong to
$PS_{D',n}$ when $n$ is large enough.  On the other hand, the limit
$u(\tsing)$ lies in $PC_{D''}$ for some finite set $D''$ strictly
contained in $D'$.  Therefore, we can apply the inductive assumption
to the sequence of ``initial data'' (\ref{main:us-sn+tn}).  We obtain
that
\begin{equation}
	\lim_{n\to +\infty}\;
	\sup_{t\geq 0}\|u_{n}(S_{n}+T_{n}+t)-u(\tsing+t)\|_{2}=0.
	\label{main:dopo-t1}
\end{equation}

Now let us observe that
\begin{eqnarray*}
	\sup_{t\geq S_{n}+T_{n}}\|u_{n}(t)-u(t)\|_{2} & \leq &
	\sup_{t\geq 0}\|u_{n}(S_{n}+T_{n}+t)-u(\tsing+t)\|_{2} \\
	 &  & + \sup_{t\geq 0}\|u(\tsing+t)-u(S_{n}+T_{n}+t)\|_{2}.
\end{eqnarray*}

The first term in the right-hand side tends to 0 as $n\to +\infty$ 
because of (\ref{main:dopo-t1}). The second term tends to 0 because 
$S_{n}+T_{n}\to\tsing$ and $u$ is uniformly continuous. It follows that
\begin{equation}
	\lim_{n\to +\infty}\;
	\sup_{t\geq S_{n}+T_{n}}\|u_{n}(t)-u(t)\|_{2}=0.
	\label{main:tn-}
\end{equation}

Therefore, (\ref{th:main}) easily follows from (\ref{main:0-sn}), 
(\ref{main:sn-tn}), (\ref{main:tn-}).

\subparagraph{\textmd{\emph{Global-in-time ``uniform'' convergence}}}

Let $w_{n}(x,t)$ be defined by
$$w_{n}(x,t):=u\left(\frac{\lfloor nx\rfloor}{n},t\right)
\quad\quad
\forall (x,t)\in[0,1]\times[0,+\infty),\ \forall n\geq 1.$$

Since $u(t)$ is already piecewise constant with respect to the space
variable, we have that $w_{n}$ and $u$ coincide in $\mathcal{K}_{n}$.
Moreover, both $u$ and $w_{n}$ are equi-bounded in $L^{\infty}$-norm,
and for every $t\geq 0$ the set $D_{n}(t)$, where $u$ and $w_{n}$ 
could be different, is the union of at most $k'$
subintervals of length $1/n$.  It follows that there exists a constant
$c_{0}$ such that
\begin{equation}
	\|u(t)-w_{n}(t)\|_{2}\leq\frac{c_{0}}{\sqrt{n}}
	\quad\quad
	\forall t\geq 0.
	\label{main:u-wn}
\end{equation}

Now for every $t\geq 0$ we have that
$$\|u_{n}(t)-w_{n}(t)\|_{2}\leq
\|u_{n}(t)-u(t)\|_{2}+
\|u(t)-w_{n}(t)\|_{2}.$$

Thanks to (\ref{th:main}) and (\ref{main:u-wn}), both terms in the
right-hand side tend to zero as $n\to +\infty$ independently of $t$,
hence
\begin{equation}
	\lim_{n\to +\infty}\;\sup_{t\geq 0}
	\|u_{n}(t)-w_{n}(t)\|_{2}=0.
	\label{main:lim-un-wn}
\end{equation}

In particular, when $n$ is large enough we have that
$$3\|u_{n}(t)-w_{n}(t)\|_{2}^{2/3}<K_{0}
\quad\quad
\forall t\geq 0,$$
where $K_{0}$ is the constant introduced in Lemma~\ref{lemma:jh-est}.
Now we exploit once again that $u$ and $w_{n}$ coincide in 
$\mathcal{K}_{n}$, and we apply (\ref{th:est-i-2}) with $v=u_{n}(t)$ 
and $w=w_{n}(t)$ (their discontinuity sets depend on time, but what 
is important is that they lie inside a fixed finite set $D$). When 
$n$ in large enough we deduce that
\begin{eqnarray*}
	\|u_{n}(x,t)-u(x,t)\|_{L^{\infty}(\mathcal{K}_{n})} & = & 
	 \sup_{t\geq 0}
	\|u_{n}(x,t)-w_{n}(x,t)\|_{L^{\infty}((0,1)\setminus D_{n}(t))}  \\
	 & \leq & \sup_{t\geq 0}
	 \|u_{n}(x,t)-w_{n}(x,t)\|_{L^{\infty}((0,1))}  \\
	 & \leq & 3\sup_{t\geq 0}
	 \|u_{n}(x,t)-w_{n}(x,t)\|_{\LL((0,1))}^{2/3}.
\end{eqnarray*}

Therefore (\ref{th:main-unif}) follows from
(\ref{main:lim-un-wn}).\qed

\subsubsection*{Proof of Theorem~\ref{cor:main}}

Let us consider the double index sequence $A_{m,n}:=u_{n}(m/n)$ with 
values in $\LL$. Due to (\ref{hp:cormain-lim-vn}) and 
(\ref{hp:cormain-lim-v}) we have that
$$A_{m,n}=v_{n}(m)\stackrel{n\to +\infty}{\longrightarrow}v(m)
\stackrel{m\to +\infty}{\longrightarrow}v_{\infty}.$$

Let us apply Lemma~\ref{lemma:amn} (refined conclusion) with
$r_{k}:=\sqrt{k}$, and let $m_{k}$ be a sequence such that
(\ref{th:amn-limit}) holds true.  Let us set $T_{n}:=m_{n}/n$.  Then
we have that $T_{n}\to 0$ and $u_{n}(T_{n})=A_{m_{n},n}\to
v_{\infty}$.

For every $n$ large enough we have that the function $u_{n}(T_{n})$
lies in some space $PS_{D''(n),n}$, with $D'\subseteq D''(n)\subseteq
D$.  Since the number of these subsets is finite, then up to
subsequences we can always assume that $D''(n)=:D''$ does not depend on
$n$.  At this point we can apply Theorem~\ref{thm:main}, with $D''$
instead of $D$, to the sequence of ``initial data'' $u_{n}(T_{n})\to
v_{\infty}$.  We obtain that
\begin{equation}
	\lim_{n\to +\infty}\,
	\sup_{t\geq 0}\|u_{n}(T_{n}+t)-u(t)\|_{2}=0.
	\label{cor:inductive}
\end{equation}

Now let us fix any $T>0$. Then for every $n$ large enough we have 
that $T_{n}\leq T$, hence
\begin{eqnarray*}
	\sup_{t\geq T}\|u_{n}(t)-u(t)\|_{2} & \leq & 
	\sup_{t\geq T_{n}}\|u_{n}(t)-u(t)\|_{2}\\
	 & = & \sup_{t\geq 0}\|u_{n}(T_{n}+t)-u(T_{n}+t)\|_{2}  \\
	 & \leq & \sup_{t\geq 0}\|u_{n}(T_{n}+t)-u(t)\|_{2}+
	 \sup_{t\geq 0}\|u(t)-u(T_{n}+t)\|_{2}.
\end{eqnarray*}

The first term in the right-hand side tends to 0 as $n\to +\infty$
because of (\ref{cor:inductive}).  The second term tends to 0 because
$u$ is uniformly continuous.  This proves (\ref{th:main-cor}).

Finally, (\ref{th:main-unif-cor}) follows from (\ref{th:main-cor}) as
(\ref{th:main-unif}) follows from (\ref{th:main}).\qed

\label{NumeroPagine}

\end{document}